\pdfoutput=1
\documentclass[12pt, oneside, reqno]{amsart}
\usepackage{amsmath,amssymb}
\usepackage{enumitem}
\usepackage{longtable}
\usepackage{booktabs,array}
\usepackage{microtype}
\usepackage{placeins}
\usepackage{needspace}
\usepackage{multicol}
\usepackage{graphicx,xcolor}
\usepackage{etoolbox}
\usepackage[textwidth=6in, textheight=9in,marginpar=0.75in]{geometry}
\usepackage[colorlinks=true,linkcolor=blue,citecolor=blue,urlcolor=blue]{hyperref}

\newtheorem{theorem}{Theorem}[section]
\newtheorem{proposition}[theorem]{Proposition}

\newtheorem{corollary}[theorem]{Corollary}

\newtheorem{maintheorem}{Theorem}

\theoremstyle{definition}

\newtheorem{remark}[theorem]{Remark}

\newcommand{\Z}{\mathbb{Z}}
\newcommand{\Q}{\mathbb{Q}}

\newcommand{\F}{\mathbb{F}}

\newcommand{\gfour}{\gamma_4}
\newcommand{\gtop}{\gamma_4^{\mathrm{top}}}
\newcommand{\ws}{\operatorname{ws}}
\newcommand{\leg}[2]{\left(\frac{#1}{#2}\right)}
\newcommand{\KI}{\textup{\textsf{KnotInfo}}}
\makeatletter
\DeclareRobustCommand{\kn}[2]{\in@{a}{#1}\ifin@\textup{#1#2}\else\in@{n}{#1}\ifin@\textup{#1#2}\else
  \ensuremath{#1_{#2}}\fi\fi}
\makeatother
\makeatletter\g@addto@macro\UrlBreaks{\do\_\do\.\do\/}\makeatother
\newcommand{\para}[1]{\par\medskip\noindent{\bfseries\boldmath#1}\par\nobreak\smallskip\noindent\ignorespaces}
\DeclareMathOperator{\Tor}{Tor}
\DeclareMathOperator{\coker}{coker}
\DeclareMathOperator{\sign}{sign}

\DeclareMathOperator{\Arf}{Arf}
\newcommand{\coderepo}{\url{https://github.com/Pinocchio315/linking-forms-nonorientable-genus}}
\newcommand{\nKnotsThirteen}{12{,}965}
\newcommand{\nHthree}{49}
\newcommand{\nMaxDetThirteen}{663}
\newcommand{\nOpen}{3{,}525}
\newcommand{\nOpenOneTwo}{991}
\newcommand{\nOpenOneThree}{773}
\newcommand{\nOpenTwoThree}{1{,}751}
\newcommand{\nOpenTwoFour}{6}
\newcommand{\nOpenOneFour}{3}

\newcommand{\nRankOneTargets}{1{,}768}
\newcommand{\nRankOneObstr}{251}
\newcommand{\nRankOneObstrOneTwo}{116}
\newcommand{\nRankOneObstrOneThree}{135}
\newcommand{\nDAltTargets}{887}
\newcommand{\nDAltObstr}{219}
\newcommand{\nDAltObstrOneTwo}{72}
\newcommand{\nDAltObstrOneThree}{147}
\newcommand{\nDMontTargets}{20}
\newcommand{\nDMontObstr}{14}
\newcommand{\nDMontObstrOneTwo}{9}
\newcommand{\nDMontObstrOneThree}{5}

\newcommand{\nNewTwoThree}{287}

\newcommand{\nThreeFromTwoThree}{41}
\newcommand{\nThreeFromOneThree}{4}

\newcommand{\nDObstrOneTwo}{81}
\newcommand{\nDObstrOneThree}{152}

\newcommand{\nOpenOneFourFive}{4}
\newcommand{\nKleinCertified}{283}
\newcommand{\nKleinUnknotBand}{194}
\newcommand{\nKleinRibbonBand}{89}
\newcommand{\nKleinTZero}{277}
\newcommand{\nKleinControlDiagrams}{9{,}312}
\newcommand{\nRangesTargets}{1{,}715}
\newcommand{\nRangesTwoThree}{1{,}710}
\newcommand{\nRangesTwoFour}{5}
\newcommand{\nRangesKI}{1{,}698}
\newcommand{\nRangesStageB}{17}
\newcommand{\nRangesTZero}{1{,}675}
\newcommand{\nRangesUnknotBand}{965}
\newcommand{\nRangesRibbonBand}{750}
\newcommand{\nElevenControlDiagrams}{226}

\newcommand{\nMobiusAll}{905}
\newcommand{\nMobiusKI}{755}
\newcommand{\nMobiusStageB}{150}
\newcommand{\nMobiusUnknot}{1}

\newcommand{\nMobiusKlein}{122}
\newcommand{\nMobiusKleinKI}{121}
\newcommand{\nMobiusControls}{484}
\newcommand{\nMobiusControlDiagrams}{91{,}026}
\newcommand{\nMobiusHOne}{879}
\newcommand{\nMobiusSig}{5}
\newcommand{\nMobiusKIgfour}{21}
\newcommand{\nMobiusFoxMilnor}{5}
\newcommand{\nMobiusFoxMilnorHolds}{16}

\newcommand{\nGammaTopTwo}{1{,}962}
\newcommand{\nGammaThreeAll}{46}
\newcommand{\nFirstOpenFromOneThree}{122}

\newcommand{\nBandsOneThree}{769}
\newcommand{\nNDTargets}{316}
\newcommand{\nNDFollowup}{62}

\newcommand{\nNDHyp}{176}
\newcommand{\nNDMixed}{30}

\newcommand{\nNDNegDef}{59}
\newcommand{\nNDNoDefAll}{51}
\newcommand{\nNDBadTwo}{9}
\newcommand{\nNDPlumbObstr}{57}

\newcommand{\nNDIndependent}{57}
\newcommand{\nNDStarRoute}{38}
\newcommand{\nNDControlsFollowup}{180}
\newcommand{\nNDControlsOne}{114}
\newcommand{\nNDControlsPipeline}{2{,}874}
\newcommand{\nNDControlsTrees}{569}
\newcommand{\nNDControlsTreeKnots}{335}
\newcommand{\nNDControlsKnown}{1{,}619}
\newcommand{\nNDControlsKnownTested}{1{,}424}
\newcommand{\nSurgTargets}{50}
\newcommand{\nSurgBox}{13}

\newcommand{\nSurgCalib}{32}
\newcommand{\nSurgCalibAgree}{31}
\newcommand{\nSurgCalibBox}{17}
\newcommand{\nSurgExactTargets}{33}

\newcommand{\nSurgTensor}{9}
\newcommand{\nSurgExactPass}{31}
\newcommand{\nSurgExactControls}{17}
\newcommand{\nSurgExcl}{19}
\newcommand{\nSurgExclOnly}{13}
\newcommand{\nNewD}{70}
\newcommand{\nNDOneTwo}{40}
\newcommand{\nNDOneThree}{30}
\newcommand{\nNDLeft}{246}
\newcommand{\nNDLeftNoMethod}{213}
\newcommand{\nSurgLeft}{31}
\newcommand{\nWideTargets}{305}
\newcommand{\nWideMobius}{12}
\newcommand{\nWideMobiusAlt}{8}
\newcommand{\nWideKI}{6}
\newcommand{\nWideReached}{10}
\newcommand{\nWideViaX}{2}
\newcommand{\nWideOpenMin}{165}
\newcommand{\nWideOpenMax}{279}
\newcommand{\nWideOpenNew}{5{,}960}
\newcommand{\nWideControls}{32}
\newcommand{\nWideControlDiagrams}{2{,}560}

\newcommand{\nDefKnots}{242}
\newcommand{\nDefSearched}{231}
\newcommand{\nDefSZero}{11}

\newcommand{\nDefDiagrams}{113{,}432}
\newcommand{\nDefWithSurface}{167}
\newcommand{\nDefCert}{64}
\newcommand{\nDefOneTwo}{35}
\newcommand{\nDefOneThree}{29}

\newcommand{\nDefCertMinC}{14}
\newcommand{\nDefCertMaxC}{23}
\newcommand{\nDefControls}{229}
\newcommand{\nDefControlsSurface}{115}
\newcommand{\nDefControlLattices}{216}
\newcommand{\nDrillTargets}{143}
\newcommand{\nDrillDescr}{21}
\newcommand{\nDrillPass}{16}
\newcommand{\nDrillDetermined}{19}
\newcommand{\nDrillNoHFK}{2}
\newcommand{\nWideAll}{16}
\newcommand{\nWideFromOneThree}{8}
\newcommand{\nKleinFromOneThreeLeft}{114}

\newcommand{\nWideFMHolds}{4}
\newcommand{\nGammaOneAll}{925}
\newcommand{\nGammaTopOneAll}{905}
\newcommand{\nMobiusOneTwoAll}{549}
\newcommand{\nMobiusOneThreeAll}{372}
\newcommand{\nGammaTwoAll}{2{,}329}
\newcommand{\nDetermined}{3{,}300}
\newcommand{\nOpenLeftOneTwo}{170}
\newcommand{\nOpenLeftFromOneThree}{55}
\newcommand{\nOpenLeft}{225}
\newcommand{\nGammaTwoOneTwo}{272}
\newcommand{\nNewDAll}{134}
\newcommand{\nKleinTotal}{2{,}120}
\newcommand{\nGammaTwoOneThree}{342}
\newcommand{\nJControls}{10}
\newcommand{\nJStopped}{16}
\newcommand{\nOpenLeftHOne}{224}

\newcommand{\nWideHOne}{9}
\newcommand{\nWideKIgfour}{5}
\newcommand{\nSharpTwo}{4{,}916}
\newcommand{\nSharpThree}{49}
\newcommand{\nSharpAll}{4{,}965}
\newcommand{\nControlsOne}{6{,}479}
\newcommand{\nDAltControls}{3{,}030}
\newcommand{\nDMontControls}{193}
\newcommand{\nValidate}{7{,}893}
\newcommand{\nGoeritzChecked}{8{,}247}
\newcommand{\nRankTwoControls}{10{,}430}
\newcommand{\nJKTargets}{1{,}028}
\newcommand{\nJKObstr}{360}
\newcommand{\nJKLinking}{141}
\newcommand{\nJKControls}{3{,}030}
\newcommand{\nJKReversed}{1{,}769}
\newcommand{\nLS}{46{,}972}
\newcommand{\nLSHthree}{196}
\newcommand{\nLSBand}{18{,}026}
\newcommand{\nLSArf}{11{,}415}
\newcommand{\nLSLinking}{15{,}431}
\newcommand{\nLSOpen}{11{,}247}
\newcommand{\nLSOpenLinking}{1{,}082}
\newcommand{\nLSDTargets}{5{,}189}
\newcommand{\nLSDObstr}{1{,}020}
\newcommand{\nLSNew}{2{,}102}
\newcommand{\nLSDControls}{5{,}987}
\newcommand{\nLSHthreeAlt}{83}
\newcommand{\nMaxDetFourteen}{1{,}145}
\newcommand{\nLSJKTargets}{5{,}754}
\newcommand{\nLSJKObstr}{1{,}585}
\newcommand{\nLSJKLinking}{565}
\newcommand{\nGapCandidates}{81}
\newcommand{\nGapAlt}{72}
\newcommand{\nGapMont}{9}
\newcommand{\nGapDiagrams}{24{,}970}
\newcommand{\nGapBands}{39{,}914{,}728}
\newcommand{\nGapKnots}{19{,}761{,}630}
\newcommand{\nGapDetOne}{22{,}453}
\newcommand{\nGapHits}{6{,}309}

\newcommand{\nGapNoForm}{8}
\newcommand{\nGapRest}{61}
\newcommand{\nGapRestNoDetOne}{6}
\newcommand{\nGapRestDetOne}{55}

\newcommand{\nGapNewCandidates}{70}
\newcommand{\nGapNew}{29}
\newcommand{\nGapNewPlumb}{19}
\newcommand{\nGapNewSurg}{10}
\newcommand{\nGapAll}{41}
\newcommand{\nGapNewKI}{13}

\newcommand{\nGapTauZeroJ}{4}

\newcommand{\nWsTwo}{17}
\newcommand{\nWsThree}{32}
\newcommand{\nWsThreeOne}{26}
\newcommand{\nWsThreeIII}{6}
\newcommand{\nLSWsTwo}{50}
\newcommand{\nMontCheckControls}{157}
\makeatletter
\def\LT@makecaption#1#2#3{%
  \LT@mcol\LT@cols c{\hbox to\z@{\hss\parbox[t]\textwidth{\reset@font\@tempcnta=64\relax
    \@makecaption{#1{#2}}{#3}}\hss}}}
\makeatother

\makeatletter
\patchcmd{\@settitle}{\bfseries}{\normalfont\scshape}{}%
  {\ClassError{amsart}{Could not set the title in small caps}{}}
\patchcmd{\@settitle}{\uppercasenonmath\@title}{}{}%
  {\ClassError{amsart}{Could not preserve title capitalization}{}}
\patchcmd{\@setauthors}{\MakeUppercase{\authors}}{{\scshape\authors}}{}%
  {\ClassError{amsart}{Could not set the author name in small caps}{}}
\g@addto@macro\@setauthors{%
  \vspace{-1.5ex}
  \begin{center}
    \normalfont\footnotesize
    Center for Geometry and Physics, Institute for Basic Science,\\
    Pohang 37673, South Korea\\
    \href{mailto:seongjinlee@ibs.re.kr}{\nolinkurl{seongjinlee@ibs.re.kr}}
  \end{center}
}
\makeatother

\hypersetup{pdftitle={Smooth and Topological Nonorientable Four-Genera of Low-Crossing Knots},pdfauthor={Seong-Jin Lee}}

\begin{document}

\title[Smooth and Topological Nonorientable Four-Genera of Low-Crossing Knots]{Smooth and Topological
Nonorientable\\ Four-Genera of Low-Crossing Knots}
\author[Seong-Jin Lee]{Seong-Jin Lee}

\begin{abstract}
We compare linking-form bounds with the smooth and locally flat nonorientable $4$-genera of low-crossing knots.  For every
prime knot with at most $13$ crossings whose Witt-level linking-form bound is at least $2$, this bound equals both
genera.  We also exhibit $\nGapAll$ prime knots with $12$ or $13$ crossings that bound locally flat M\"obius bands but
no smooth ones, and whose smooth nonorientable $4$-genus is $2$.  Our methods combine linking forms, Heegaard Floer
correction terms and definite-lattice obstructions with explicit band constructions.  They determine $\nDetermined$ of
the $\nOpen$ values left unknown in \KI{} on 2 November 2026; the remaining $\nOpenLeft$ knots have range $[1,2]$.
We obtain new lower bounds of $2$ for $\nLSNew$ knots with $14$ crossings and correct previously reported values for
three knots with $11$ crossings.
\end{abstract}


\maketitle
\section{Introduction}\label{sec:intro}

The \emph{nonorientable $4$-genus} $\gfour(K)$ of a knot $K\subset S^3$ is the least first Betti number of a properly
and smoothly embedded compact connected nonorientable surface in $B^4$ with boundary $K$, with the convention that
$\gfour(K)=0$ if $K$ is smoothly slice \cite{MY}.  Its locally flat analogue $\gtop(K)$ is defined in the same way, with
value $0$ for topologically slice knots.  Thus $\gtop(K)\le\gfour(K)$.  We study how linking forms and smooth
obstructions determine these genera when combined with explicit surfaces for knots of small crossing number.
Let $h_\lambda(K)$ be the Witt-level linking-form bound of \cite{LeeTheory}, recalled in \S\ref{sec:hlambda}.  It
replaces the orthogonal-summand criterion of Murakami and Yasuhara \cite{MY} and Gilmer and Livingston \cite{GL}, which
is false in general \cite{LeeTheory}, and for the knots with at most $14$ crossings it keeps every obstruction drawn
from that criterion (\S\ref{sec:tools}).  It satisfies
\[
h_\lambda(K)\le\gtop(K)\le\gfour(K).
\]
Our first result shows that both inequalities are equalities for every prime knot with at most $13$ crossings and
$h_\lambda(K)\ge2$.

\begin{maintheorem}\label{mthm:sharpness}
For every prime knot $K$ with at most $13$ crossings and $h_\lambda(K)\ge2$,
\[
\gtop(K)=\gfour(K)=h_\lambda(K).
\]
There are $\nSharpAll$ such knots: $\nSharpTwo$ with $h_\lambda=2$ and $\nSharpThree$ with $h_\lambda=3$.
\end{maintheorem}

This statement includes knots whose genera were previously known.  Its proof combines the linking-form lower bound with
the upper bounds established here and those recorded in \KI{} \cite{KnotInfo}; see \S\ref{sec:values}, including the
corrections in \S\ref{sec:corrections}.  All lower bounds of $3$ used here follow from the linking form, signature and
Arf invariant through Corollary~\ref{cor:criteria}(c).

At $h_\lambda=1$, smooth obstructions detect a difference between the two categories.
\begin{maintheorem}\label{mthm:gap}
Each of the $\nGapAll$ prime knots of Table~\ref{tab:gap}, all with $12$ or $13$ crossings, bounds a locally flat
M\"obius band in $B^4$ but no smooth one.  For each of them,
\[
h_\lambda(K)=\gtop(K)=1<\gfour(K)=2.
\]
\end{maintheorem}

Gilmer and Livingston constructed knots bounding locally flat but no smooth M\"obius bands \cite[Thm.~D]{GL}, and Feller
and Golla established $\gtop<\gfour$ for the torus knots $T_{2n-1,2n}$ with $n\ge5$ \cite[Thm.~1.3]{FG}.  Prime examples
also follow from Gilmer and Livingston's result and Myers' theorem that every nontrivial knot is smoothly concordant to
a hyperbolic knot \cite[Thm.~1.1]{Myers}, since both genera are invariant under smooth concordance.
Theorem~\ref{mthm:gap} gives explicit examples with $12$ and $13$ crossings and determines both genera.  For each
example, a band move to a knot with Alexander polynomial one gives the locally flat upper bound, while correction terms
exclude a smooth M\"obius band; a smooth surface with first Betti number $2$ completes the determination
(\S\ref{sec:gap}).

These results form part of a systematic computation of the nonorientable $4$-genus through $13$ crossings.  Earlier
computations include those of Jabuka and Kelly \cite{JK}, Ghanbarian \cite{Ghanbarian}, Fairchild \cite{Fairchild}, and
Lee and Sabloff \cite{LS}.  We use the values and ranges recorded in \KI{}, version 2026.11.2, which agree on the data
used here with the earlier versions used for the computations (\S\ref{sec:methods}).  \KI{} assigns value $1$ to slice
knots; our convention assigns value $0$.

\begin{maintheorem}\label{mthm:main}
Of the $\nOpen$ prime knots with at most $13$ crossings for which \KI{} 2026.11.2 lists a range, we determine $\gfour$
for $\nDetermined$: it is $1$ for $\nGammaOneAll$ knots, $2$ for $\nGammaTwoAll$, and $3$ for $\nGammaThreeAll$.  Among
these, $\gtop=1$ for $\nGammaTopOneAll$ of the knots with $\gfour=1$, and $\gtop=2$ for $\nGammaTopTwo$ and $\gtop=1$
for $\nGapAll$ of the knots with $\gfour=2$.  Every knot with $\gfour=3$ in this determination has $\gtop=3$.  Every
range with lower bound $2$ is settled.  The remaining $\nOpenLeft$ knots all have range $[1,2]$;
$\nOpenLeftFromOneThree$ of them previously had range $[1,3]$.
\end{maintheorem}

The lower and upper bounds are obtained by separate methods.  For smooth M\"obius bands, Proposition~\ref{prop:dtest}
gives a correction-term obstruction that accounts for isotropic subgroups of the linking form.  We obtain correction
terms, or bounds sufficient for this test, from negative-definite plumbings and surgery descriptions of the double
branched cover.  When correction terms are unavailable, definite checkerboard surfaces in non-minimal diagrams allow us
to apply Donaldson's theorem (Proposition~\ref{prop:donaldson}).  Together these methods give $\nNewDAll$ further values
for non-alternating knots that are not Montesinos (Proposition~\ref{prop:newd}).  We do not use the bound
$|\upsilon(K)-\sigma(K)/2|\le\gfour(K)$ of Ozsv\'ath, Stipsicz and Szab\'o \cite[Thm.~1.2]{OSS}, where $\upsilon$ is a
concordance invariant from knot Floer homology and $\sigma$ is the signature: it gives nothing for quasi-alternating
knots, for which $\upsilon(K)=\sigma(K)/2$ \cite[p.~3]{OSS}, and we have not computed $\upsilon$ for the other knots.

On alternating knots, the relation between the correction-term and lattice tests is explained by sharp fillings and
Elkies' theorem (\S\ref{sec:donaldson}).  The new upper bounds come from band constructions: the resulting sequences of
band and Reidemeister moves end in unlinks, and a separate computation rebuilds the surfaces.  In particular, the
constructed surfaces do not depend on the obstruction calculations or on previously recorded values of $\gfour$
(\S\ref{sec:upper}).

The lower-bound computations also extend to $14$ crossings.
\begin{maintheorem}\label{mthm:fourteen}
Of the $\nLSOpen$ knots with $14$ crossings left open by \cite{LS}, $\nLSNew$ have $\gfour\ge2$.  The $\nLSHthree$ knots
with $14$ crossings and $h_\lambda=3$ in the table of \cite{LS} have $\gtop\ge3$.
\end{maintheorem}

The computations also settle or clarify several values at lower crossing number.
\begin{maintheorem}\label{mthm:corrections}
We have
\[
\gfour(\kn{11n}{133})=\gtop(\kn{11n}{133})=3,\qquad\gfour(\kn{11n}{57})=2.
\]
The published argument for the lower bound $2$ for $\kn{11n}{176}$ fails, and its value remains in $[1,2]$.  For
$\kn{11n}{38}$, the value $2$ holds with a different proof of the lower bound.  We also have
$\gfour(\kn{9}{49})=\gtop(\kn{9}{49})=3$.
\end{maintheorem}

The value $\gfour(\kn{9}{49})=3$ was already established in the revised version of Jabuka and Kelly \cite{JKv3}.  For
$\kn{11n}{57}$, Lee and Sabloff had already found that the band move of \cite[Prop.~3.2]{Fairchild} does not give the
unknot \cite[\S4.3.2]{LS}; our lower bound completes its determination.  The individual arguments and their relation to
the published proofs are given in \S\ref{sec:corrections}.

Section~\ref{sec:tools} gives the bounds and the band conventions.  Section~\ref{sec:computations} proves
Theorems~\ref{mthm:sharpness} and~\ref{mthm:main}, explains the checks, and treats the knots with $14$ crossings and the
corrections (Theorems~\ref{mthm:fourteen} and~\ref{mthm:corrections}).  Section~\ref{sec:gap} gives the locally flat
M\"obius bands of Theorem~\ref{mthm:gap}.  Appendix~\ref{app:tables} contains the knot lists.  The accompanying code and
data record the diagrams, calculations and verification results used in the proofs.

\section{The bounds}\label{sec:tools}

\subsection{The linking-form bound}\label{sec:hlambda}
For a knot $K$ let $\Sigma_K$ be its double branched cover, $H_K=H_1(\Sigma_K)$ and $D=\det K=|H_K|$, which is odd.  A
nondegenerate symmetric matrix $Q\in M_h(\Z)$ has the \emph{discriminant form}
\[
\lambda_Q(\bar x,\bar y)=x^{T}Q^{-1}y\bmod\Z\qquad\text{on }\coker Q=\Z^h/Q\Z^h,
\]
where $\bar x$ is the class of $x\in\Z^h$.  For a Seifert matrix $V$ of $K$ we put $\lambda_K=\lambda_{V+V^{T}}$, a
linking form on a group isomorphic to $H_K$; with the sign convention \cite[(6)]{LeeTheory} it is minus
the linking form of $\Sigma_K$ \cite[\S2.3]{LeeTheory}.  Let $\sigma(K)=\sign(V+V^{T})$ be the
signature of $K$ and $\Arf(K)$ its Arf invariant.  For an odd prime $p$ let $\lambda_{K,p}$ be the $p$-primary part of
$\lambda_K$, let $a_p\in\{0,1,2\}$ and $\delta_p\in\F_p^\times/\F_p^{\times2}$ be the rank and the discriminant (the class of the
determinant, with no sign factor) of the anisotropic quadratic space over $\F_p$ that is Witt equivalent to the odd-level
part of $\lambda_{K,p}$ \cite[\S2.1]{LeeTheory}, and let
$D'_p$ be the prime-to-$p$ part of $D$.  For $h\ge1$ let $\mathcal Q_h$ be the set of odd symmetric $Q\in M_h(\Z)$ of
odd determinant, $Q$ being \emph{odd} if some diagonal entry is odd, and let $\mathcal Q_0$ consist of the empty
matrix, for which $\lambda_Q=0$.  Writing $[\lambda]$ for the class of $\lambda$ in the Witt group of linking forms on
groups of odd order, we put
\[
h_\lambda(K)=\min\bigl\{h\ge0:\ [\lambda_Q]=[\lambda_K]\text{ for some }Q\in\mathcal Q_h\bigr\}
\]
\cite[Def.~4.1]{LeeTheory}.  We use the following results of \cite{LeeTheory}.

\begin{theorem}[{\cite[Thm.~4.2]{LeeTheory}}]\label{thm:lower}
For every knot $K$, $h_\lambda(K)\le\gtop(K)\le\gfour(K)$.
\end{theorem}

\begin{theorem}[{\cite[Thm.~4.4]{LeeTheory}}]\label{thm:compute}
Let $h\ge0$.  Then $h_\lambda(K)\le h$ if and only if there is $\tau\in\{\pm1\}$ such that the following two conditions
hold.
\begin{enumerate}[label=(L$_h$),leftmargin=3.2em]
\item $a_p\le h$ for every odd prime $p$, and $\leg{\delta_p}{p}=\leg{\tau D'_p}{p}$ for every $p$ with $a_p=h\ge1$.
\end{enumerate}
\begin{enumerate}[label=(R$_h$),leftmargin=3.2em]
\item $h\ne2$; or $\tau=1$; or $\sigma(K)+4\Arf(K)\not\equiv4\pmod8$.
\end{enumerate}
In particular $h_\lambda(K)\le3$ for every knot.
\end{theorem}

\begin{corollary}[{\cite[Cor.~4.5]{LeeTheory}}]\label{cor:criteria}
For a prime $p$ with $a_p\ge1$ let $\varepsilon_p=\leg{\delta_pD'_p}{p}$.
\begin{enumerate}[label=(\alph*)]
\item $h_\lambda(K)\ge1$ if and only if $\lambda_K$ is not metabolic.
\item $h_\lambda(K)\ge2$ if and only if either some $a_p=2$, or for each $\tau\in\{\pm1\}$ some $p$ with $a_p=1$ has
$\leg{\tau}{p}\ne\varepsilon_p$.  In particular $h_\lambda(K)\ge2$ whenever $\sigma(K)+4\Arf(K)\equiv4\pmod8$.
\item \emph{(Klein bottle criterion.)} $h_\lambda(K)\ge3$ if and only if one of the following holds:
  \begin{enumerate}[label=(\roman*)]
  \item some $p\equiv1\pmod 4$ with $a_p=2$ has $\varepsilon_p=-1$;
  \item two primes $p,q\equiv3\pmod4$ with $a_p=a_q=2$ have $\varepsilon_p\ne\varepsilon_q$;
  \item some $p\equiv3\pmod4$ with $a_p=2$ has $\varepsilon_p=-1$, and $\sigma(K)+4\Arf(K)\equiv4\pmod8$.
  \end{enumerate}
\end{enumerate}
Consequently, if one of (i)--(iii) holds, then $\gtop(K)\ge3$: $K$ bounds no locally flat M\"obius band and no locally
flat punctured Klein bottle in $B^4$.
\end{corollary}

The last statement of Corollary~\ref{cor:criteria}(b) is Yasuhara's condition \cite[Prop.~5.1]{Yasuhara}.  The splitting
lemma would give instead an orthogonal summand presented by an odd form \cite[\S4.2]{LeeTheory}: let
$h^\oplus(K)$ be the least $h\ge0$ such that $\lambda_K\cong\lambda_Q\oplus\mu$ for some $Q\in\mathcal Q_h$ and a
metabolic $\mu$ \cite[Def.~4.8]{LeeTheory}.  A prime $p$ is
\emph{exceptional} for $\lambda_K$ if $p\equiv3\pmod 4$ and the levels of odd exponent in a Jordan decomposition of
$\lambda_{K,p}$ \cite[\S2.1]{LeeTheory} all have rank one, have discriminants in one square class and are
$\equiv3\pmod 4$ in number \cite[Def.~3.11]{LeeTheory}; then $p^9\mid D$.

\begin{proposition}[{\cite[Prop.~4.9]{LeeTheory}}]\label{prop:hoplus}
If no prime is exceptional for $\lambda_K$, then $h^\oplus(K)=h_\lambda(K)$.
\end{proposition}

The knots with at most $14$ crossings have determinant at most $\nMaxDetFourteen<3^9$, so no prime is exceptional for
them and $h^\oplus=h_\lambda$; in particular the obstructions drawn from the criterion of \cite{MY,GL} remain valid for
them.  The same holds for all knots with at most $16$ crossings (Remark~\ref{rem:sixteen}).  The criterion itself allows
forms of either parity and is weaker: $\lambda_{\kn{8}{18}}$ is presented by an even binary form, while
$h_\lambda(\kn{8}{18})=3$ \cite[Example~4.7]{LeeTheory}.  In general $h^\oplus$ and $h_\lambda$ differ, and
$h^\oplus$ is not a lower bound for $\gtop$ \cite[Thm.~5.1]{LeeTheory}.

\begin{remark}\label{rem:sixteen}
Every knot $K$ with at most $16$ crossings has $\det K<3^9$, so no prime is exceptional for $\lambda_K$ and
$h^\oplus(K)=h_\lambda(K)$ (Proposition~\ref{prop:hoplus}).  Indeed, colour a diagram of $K$ with $c$ crossings as a
checkerboard, and let $X_0,\dots,X_{v-1}$ be the white regions and $\eta(D)=\pm1$ the incidence number of a crossing $D$.
The Goeritz matrix $G=(g_{ij})_{1\le i,j\le v-1}$ is given by
\[
g_{ij}=-\sum_{D}\eta(D)\quad(i\ne j),\qquad g_{ii}=-\sum_{k\ne i}g_{ik},
\]
where $D$ runs over the crossings incident to $X_i$ and $X_j$, and $k$ over $0,\dots,v-1$ \cite[\S1, p.~55]{GoLi};
and $|\det G|=\det K$, as $G$ presents $H_K$ \cite[p.~60]{GoLi}.  So $G$ is the Laplacian of the graph $\Gamma$ whose
vertices are the white regions and whose edges are the crossings, weighted by $\eta$, with the row and the column of
$X_0$ removed.  By Kirchhoff's matrix-tree theorem $\det G=\sum_T\prod_{D\in T}\eta(D)$, summed over the spanning
trees $T$ of $\Gamma$, and a spanning tree consists of $v-1$ of the $c$ edges; so
\[
\det K\le\#\{T\}\le\binom{c}{v-1}\le\binom{c}{\lfloor c/2\rfloor}\le\binom{16}{8}=12870<3^9 .
\]
In particular a knot with an exceptional prime has at least $17$ crossings.
\end{remark}

\para{The examples of determinant $27$.}
In \cite[Examples~3.5 and~3.6]{LeeTheory} the splitting lemma fails for a smooth
M\"obius band and a smooth punctured torus whose boundaries $K$ and $K'$ have $H_K\cong H_{K'}\cong\Z/27$; there $K$ is
obtained from the pretzel knot $P(3,-3,1)$ by the band move along a band $b_{-3}$, and $K'$ by the band moves along
$b_{-2}$ and along a band $b'$ in the bigon between the two crossings of $b_{-2}$.  We built diagrams of $K$ and $K'$
in this way.  Their PD codes, in the convention of \KI{} \cite{KnotInfo}, in which the edges are numbered from $1$ along
the orientation and each crossing lists its four edges counterclockwise, beginning with the incoming under-strand, are
\begin{align*}
K\colon\ &(20,14,1,13),\ (14,2,15,1),\ (2,16,3,15),\ (5,12,6,13),\ (11,4,12,5),\\
&(3,10,4,11),\ (16,9,17,10),\ (17,7,18,6),\ (7,19,8,18),\ (19,9,20,8);\\
K'\colon\ &(22,9,1,10),\ (10,21,11,22),\ (20,11,21,12),\ (17,9,18,8),\\
&(7,19,8,18),\ (19,7,20,6),\ (12,5,13,6),\ (13,16,14,17),\\
&(1,4,2,5),\ (2,15,3,16),\ (14,3,15,4).
\end{align*}
Their Goeritz matrices on the white regions, with the incidence numbers of Gordon and Litherland, are those of
\cite{LeeTheory}, sign included, and SnapPy \cite{SnapPy} finds isometries between their complements and those of
$\kn{8}{10}$ and $\kn{8}{11}$.  So $K$ and $K'$ are $\kn{8}{10}$ and $\kn{8}{11}$ up to mirror image
\cite[Thm.~1]{GoLu}.

\subsection{Correction terms}\label{sec:dtest}

Theorem~\ref{thm:lower} concerns locally flat surfaces.  A smooth M\"obius band satisfies a further condition, which
comes from the correction terms $d$ of Heegaard Floer homology \cite{OSz-abs}; Gilmer and Livingston used them to
obstruct smooth punctured Klein bottles \cite[\S5]{GL} and, in the proof of \cite[Thm.~D]{GL}, smooth M\"obius bands.
Since $D$ is odd, $\Sigma_K$ has a unique spin structure $\mathfrak t_s$, and every Spin$^c$ structure is $\mathfrak
t_s+x$ for a unique $x\in H^2(\Sigma_K)$, which we identify with $H_K$ by Poincar\'e duality.  We write $\lambda_Z$ for
the linking form of $Z=\pm\Sigma_K$; thus $\lambda_{\Sigma_K}=-\lambda_K$ and $\lambda_{-\Sigma_K}=\lambda_K$ (\S\ref{sec:hlambda}).

\begin{proposition}\label{prop:dtest}
Suppose that $K$ bounds a smooth M\"obius band in $B^4$.  Then for $Z=\Sigma_K$ or $Z=-\Sigma_K$ there are an odd
integer $\ell\ge1$, an isotropic subgroup $M\subset H_K$ with $\ell|M|^2=D$, and an element $g\in M^\perp$ that generates
$M^\perp/M\cong\Z/\ell$ and has $\lambda_Z(g,g)=1/\ell$, such that
\begin{equation}\label{eq:dtest}
4\,d(Z,\mathfrak t_s+jg+m)\ \ge\ 1-\frac{(2j-\ell)^2}{\ell}\qquad\text{for all integers }0\le
j<\ell\text{ and all }m\in M .
\end{equation}
\end{proposition}

\begin{proof}
Let $F$ be a smooth M\"obius band with $\partial F=K$ and $W=\Sigma_2(B^4,F)$.  By
\cite[Lemma~3.3]{LeeTheory}, $b_2(W)=1$ and $H_1(W)$ is finite of odd order, so $Q_W=(q)$ with $q$ odd;
put $\ell=|q|$.  Let $W_0=W$ if $q<0$ and $W_0=-W$ if $q>0$.  Then $W_0$ is smooth and negative definite,
$Q_{W_0}=(-\ell)$, and $\partial W_0=Z$ with $Z=\pm\Sigma_K$.  Let $e$ generate $H_2(W_0)/\Tor$, let $e^*\in H^2(W_0)$
satisfy $\langle e^*,e\rangle=1$, and let $g\in H^2(Z)=H_K$ be its restriction.  Under Lefschetz duality, restriction to
$Z$ is the boundary map of \cite[Lemma~3.1]{LeeTheory}, and the image of $\Tor H^2(W_0)$ is the
subgroup $M$ of that lemma.  So $M$ is isotropic, $\ell|M|^2=D$, and $g\in M^\perp$ maps to a generator of
$M^\perp/M\cong\coker Q_{W_0}$, with $\lambda_Z(g,g)=-\lambda_{Q_{W_0}}(1,1)=1/\ell$.

Since $H^2(W_0)$ has no $2$-torsion, a Spin$^c$ structure $\mathfrak s$ on $W_0$ is determined by $c_1(\mathfrak s)$.
The characteristic classes are $ke^*+\theta$ with $k$ odd and $\theta\in\Tor H^2(W_0)$.  For such $\mathfrak s$ we have
$c_1(\mathfrak s)^2=-k^2/\ell$ and $c_1(\mathfrak s|_Z)=kg+\theta|_Z$, so $\mathfrak s|_Z=\mathfrak
t_s+\frac12(kg+\theta|_Z)$, halving in the group $H_K$ of odd order.  Since $\ell g\in M$, these restrictions fill the
coset $\mathfrak t_s+\frac{k+\ell}{2}g+M$ as $\theta$ varies.  By the inequality of Ozsv\'ath and Szab\'o
\cite[Thm.~9.6]{OSz-abs}, restated as \cite[Thm.~12]{GL}, a smooth negative-definite $4$-manifold $X$ whose boundary
$Y$ is a rational homology sphere satisfies
\[
c_1(\mathfrak s)^2+b_2(X)\le4d(Y,\mathfrak s|_Y)\qquad\text{for every Spin$^c$ structure $\mathfrak s$ on $X$.}
\]
For $X=W_0$ this reads $c_1(\mathfrak s)^2+1\le4d(Z,\mathfrak s|_Z)$, and with $k=2j-\ell$ it is \eqref{eq:dtest}.
\end{proof}

\para{Sharp fillings.}
For an alternating knot, $\Sigma_K$ or $-\Sigma_K$ bounds the negative-definite $4$-manifold $X_G$ of a Goeritz matrix
$G$ of rank $n$, and $X_G$ is sharp \cite[Thm.~3.4]{OSz-bdc}:
\[
d(\partial X_G,\mathfrak t)=\max_c\ \frac{c^{T}G^{-1}c+n}{4},
\]
the maximum over the characteristic covectors $c\in\Z^n$, those with $c_i\equiv G_{ii}\pmod 2$ for all $i$, in the
class of $\mathfrak t$; here the Spin$^c$ structures on $\partial X_G$ correspond to the characteristic covectors
modulo $2G\Z^n$.  For a Montesinos knot, $\Sigma_K$ or $-\Sigma_K$ bounds a negative-definite star-shaped plumbing,
which has at most one bad vertex, and the same formula holds for its intersection matrix \cite[Thm.~1.2 and
Cor.~1.5]{OSz-plumb}.  More generally, it holds whenever $\Sigma_K$ or $-\Sigma_K$ bounds a negative-definite plumbing
tree with at most two bad vertices, a vertex being bad if its weight exceeds minus its degree
\cite[Cor.~1.5]{OSz-plumb}.  In these cases Proposition~\ref{prop:dtest} is a finite computation.

\para{Surgery descriptions.}
If $K$ has unknotting number one, then $\Sigma_K\cong\pm S^3_{D/2}(C)$ for a knot $C\subset S^3$, by the Montesinos trick
\cite{Montesinos}, \cite[Lemma~1.5]{OSz-unknot}.  More generally, if drilling a closed curve from $\Sigma_K$ gives the
exterior of a knot $C\subset S^3$, then $\Sigma_K\cong S^3_{p/q}(C)$ for a slope $p/q$ with $|p|=D$.
For $p,q>0$ the correction terms of $S^3_{p/q}(C)$ are
\[
d\bigl(S^3_{p/q}(C),i\bigr)=d\bigl(L(p,q),i\bigr)-2\max\bigl\{V_{\lfloor i/q\rfloor}(C),\,
V_{\lfloor(p+q-1-i)/q\rfloor}(C)\bigr\},\qquad 0\le i<p,
\]
by \cite[Prop.~1.6 and Rem.~2.10]{NiWu}; a negative slope reduces to this case, as
$S^3_{-p/q}(C)\cong-S^3_{p/q}(\overline C)$ for the mirror image $\overline C$ of $C$.  Here $L(p,q)=S^3_{p/q}(O)$.
With this orientation convention
\cite[\S1.2]{OSz-lens}, the recursion of \cite[Prop.~4.8]{OSz-abs} and the relation $d(-Y,\mathfrak t)=-d(Y,\mathfrak
t)$ \cite[Prop.~4.2]{OSz-abs} give, for $p>q>0$,
\[
d\bigl(L(p,q),i\bigr)=\frac{(2i+1-p-q)^2-pq}{4pq}-d\bigl(L(q,r),j\bigr),\qquad 0\le i<p,
\]
where $r$ and $j$ are the reductions of $p$ and $i$ modulo $q$, and $d(S^3)=0$ ends the recursion.  The integers
$V_j(C)\ge0$ are invariants of the knot Floer complex of $C$; and $i$ labels the Spin$^c$ structures as in
\cite[\S7]{OSz-rat}, where the meridian $\mu$ of $C$ corresponds to $q$ \cite[proof of Lemma~7.2]{OSz-rat}, so that
$\mathfrak t_s+k\mu$ has the label $i_s+qk$ with $2i_s\equiv q-1\pmod p$.  In the sign
convention of the proof of Proposition~\ref{prop:dtest}, in which the boundary of a $4$-manifold with intersection
matrix $A$ has the linking form $-A^{-1}$, $\lambda(\mu,\mu)=-q/p$.

\para{The integers $V_j(C)$.}
When the $V_j(C)$ are not known, we replace them by the bounds
\[
V^-_j\le V_j(C)\le V^+_j,\qquad V^-_j=\begin{cases}1,&0\le j<\tau(C),\\0,&j\ge\tau(C),\end{cases}\qquad
V^+_j=\max(g_S-j,0).
\]
The lower bound holds since $V_{j+1}(C)\le V_j(C)$ and $\nu^+(C)=\min\{j\ge0:V_j(C)=0\}\ge\tau(C)$ \cite[(2.3),
Def.~2.1 and Prop.~2.3]{HomWu}; the upper bound is \cite[Cor.~7.4]{Rasmussen}, where $g_S\ge g_4(C)$ is the genus of
the Seifert surface that Seifert's algorithm gives on a diagram of $C$.  The lower bound needs $\tau(C)$, read off the
knot Floer homology of $C$ computed with Szab\'o's program \cite{Szabo}, and the chirality of $C$.  Each choice of
$Z$, $\ell$, $M$ and $g$ in Proposition~\ref{prop:dtest} imposes inequalities $4d(Z,\mathfrak t)\ge b$, and
$d(S^3_{p/q}(C),\mathfrak t)$ does not increase when some $V_j$ increases; so the choice satisfies \eqref{eq:dtest} for
some $V$ between the two bounds if and only if it does so for $V^-$ when $Z=S^3_{p/q}(C)$, and for $V^+$ when
$Z=-S^3_{p/q}(C)$.  The two computations, with $V^-$ and with $V^+$, therefore decide whether the test excludes a
smooth M\"obius band for every $V$ between the bounds; if it does not, the bounds alone decide nothing for the actual
$V(C)$.
The exact values come from the knot Floer complex of $C$.  Szab\'o's program gives a model of it over
$\F_2[U,V]/(UV)$ with $\dim\widehat{\mathit{HFK}}(C)$ generators; as setting $U=V=0$ gives a complex with homology
$\widehat{\mathit{HFK}}(C)$, the model is reduced, that is, its differential vanishes modulo $(U,V)$.  The integer
$V_s(C)$ is read off the homology of the quotient complex
$A^+_s=\mathit{CFK}^\infty(C)\{\max(i,j-s)\ge0\}$, whose tower has its bottom in grading $-2V_s(C)$ \cite{NiWu}, and we determine it in two
ways.  First, multiplication by $U$ maps $A^+_s$ onto itself with kernel
$\widehat A_s=\mathit{CFK}^\infty(C)\{\max(i,j-s)=0\}$, so that $0\to\widehat A_s\to A^+_s\xrightarrow{U}A^+_s\to0$ is
exact and the bottom of the tower comes from $H_{-2V_s}(\widehat A_s)$; and $\widehat A_s$ contains no term $U^aV^b$ with $a,b\ge1$,
so its homology can be computed from the reduced model, and it determines $V_s(C)$ when it is nonzero in exactly one
even grading $\le0$.  Second, the complex over $\F_2[U,V]$ has a reduced model that agrees with the given one modulo
$UV$, in the same basis.  Indeed, a homotopy equivalence between reduced complexes over $\F_2[U,V]/(UV)$ is invertible
modulo $(U,V)$, where both differentials vanish, so it is an isomorphism by the graded Nakayama lemma, and a lift to
$\F_2[U,V]$ of a graded isomorphism is invertible for the same reason.  So the reduced model over $\F_2[U,V]$ differs
from the given one only by terms $U^aV^b$ with $a,b\ge1$ that both gradings allow; we compute $V_s$ for every choice of
such terms with $\partial^2=0$.

\subsection{Donaldson's theorem}\label{sec:donaldson}

Donaldson's theorem gives a condition of a different kind, through spanning surfaces whose Gordon--Litherland form is
definite.  For a diagram of $K$ coloured as a checkerboard, the Goeritz matrix of \cite[\S1]{GoLi}, formed on the white
regions with the incidence numbers $\eta$ (Remark~\ref{rem:sixteen}), is the matrix of the Gordon--Litherland form
$\mathcal G_V$ of the black surface $V$ \cite[Thm.~1]{GoLi}, and the double cover $X_V$ of $B^4$ branched along $V$
pushed into $B^4$ has intersection form $(H_1(V),\mathcal G_V)$ \cite[Thm.~3]{GoLi}.

\begin{proposition}\label{prop:donaldson}
Suppose that $K$ bounds a smooth M\"obius band in $B^4$, that $\sigma(K)+4\Arf(K)\equiv2\varepsilon\pmod 8$ with
$\varepsilon=\pm1$, and that $V$ is a spanning surface of $K$ whose Gordon--Litherland form $\mathcal G_V$ is negative
definite if $\varepsilon=1$ and positive definite if $\varepsilon=-1$.  Then the positive definite lattice
$(H_1(V),-\varepsilon\mathcal G_V)$ embeds in the standard lattice $\Z^{b_1(V)+1}$.
\end{proposition}

\begin{proof}
Let $F$ be a smooth M\"obius band with $\partial F=K$ and $W=\Sigma_2(B^4,F)$, so that $Q_W=(q)$ with $q\ne0$, as in
the proof of Proposition~\ref{prop:dtest}.  By \cite[Cor.~8]{GL},
\[
\sigma(K)+4\Arf(K)\equiv\sign W-\beta(B^4,F)\pmod 8,
\]
where $\beta$ is the Brown invariant; for a M\"obius band $\sign W=\pm1$ and $\beta(B^4,F)\equiv\pm1\pmod 8$
\cite[proof of Thm.~5]{GL}, so $\sign W=\varepsilon$, that is, $q$ has the sign of $\varepsilon$.  The boundaries of
$X_V$ and $W$ are both $\Sigma_K$, so $N=X_V\cup_{\Sigma_K}(-W)$ is a closed smooth oriented $4$-manifold.  Since
$\Sigma_K$ is a rational homology sphere, $H_2(\Sigma_K)=0$, and by the Mayer--Vietoris sequence
$(H_2(X_V)\oplus H_2(-W))/\Tor$ maps injectively onto a subgroup of finite index of $H_2(N)/\Tor$, on which the
intersection form is $\mathcal G_V\oplus(-q)$.  So $N$ is negative definite if $\varepsilon=1$ and positive definite if $\varepsilon=-1$, and by Donaldson's theorem \cite[Thm.~1]{Donaldson}, applied to $N$ or to $-N$, its intersection form is $-\varepsilon$ times the standard form of rank $b_1(V)+1$.  Its restriction to $H_2(X_V)$, which is free, is
the required embedding.
\end{proof}

For a checkerboard surface this is the criterion of Jabuka and Kelly \cite[Thm.~2.10]{JK}, which they state for every
knot with a definite Goeritz form; in their statement the diagonal form has the sign $\varepsilon$, in which no
$(-\varepsilon)$-definite form embeds, and the argument gives the sign used here.

For alternating knots the two criteria agree.  Such a knot has sharp
negative-definite Goeritz manifolds bounded by $\Sigma_K$ and by $-\Sigma_K$, from an alternating diagram of $K$ and of its
mirror image \cite[Thm.~3.4]{OSz-bdc}.  Let $X$ be the one with $\partial X=-Z$, of rank $n$,
and for $\ell$, $M$ and $g$ as in Proposition~\ref{prop:dtest} let $\Lambda$ be the unimodular lattice that contains
$H_2(X)\oplus\langle-\ell\rangle$ and is defined by the anti-isometry between $M^\perp/M$ and the discriminant group of
$\langle-\ell\rangle$ that $g$ determines.  Explicitly, the discriminant groups of $H_2(X)$ and of $\langle-\ell\rangle$
are $H_K$ with the form $\lambda_Z$ and $\Z/\ell$ with the form $-1/\ell$, and $\Lambda/(H_2(X)\oplus\langle-\ell\rangle)$ is
\[
S=\{(x,j)\in H_K\oplus\Z/\ell:\ x\in jg+M\},
\]
which is well defined since $\ell g\in M$, is isotropic since $M$ is isotropic, $g\in M^\perp$ and $\lambda_Z(g,g)=1/\ell$,
and has order $\ell|M|=\sqrt{D\ell}$, so that $\Lambda$ is unimodular.  The characteristic covectors of $\Lambda$ are the pairs $(c,ke^*)$ of
characteristic covectors of $H_2(X)$ and of $\langle-\ell\rangle$ whose classes on $Z$ correspond as in the proof of
Proposition~\ref{prop:dtest}, and such a pair has square $c^2-k^2/\ell$.  Since $X$ is sharp, $4d(Z,\mathfrak
t)=-\max(c^2+n)$ over the class of $\mathfrak t$, so \eqref{eq:dtest} holds for all $j$ and $m$ if and only if
\[
\xi^2\le-(n+1)\qquad\text{for every characteristic covector $\xi$ of $\Lambda$}.
\]
By Elkies' theorem \cite{Elkies} this is the case if and only if $\Lambda\cong\langle-1\rangle^{\oplus(n+1)}$, that is, if
$H_2(X)\oplus\langle-\ell\rangle$ embeds in the diagonal lattice of rank $n+1$ with this gluing.  If an embedding sends the
generator of $\langle-\ell\rangle$ to $s$ times a primitive vector $w$, then $w$ is orthogonal to $H_2(X)$, has square
$-\ell/s^2$, and $H_2(X)\oplus\langle w\rangle$ is glued in this way, with $\ell/s^2$ in place of $\ell$.  Hence,
for alternating knots, Proposition~\ref{prop:dtest} excludes the same knots as the lattice criterion of \cite[Thms.~2.10
and 2.12]{JK} applied with both signs of $Q_W$; Greene \cite{GreeneNote} and Larson \cite[Thm.~1.1]{Larson} use the same
argument for rational balls.
For the alternating knots with at most $14$ crossings that we consider, the criterion of \cite{JK}, with the sign of $Q_W$ fixed by
$\sigma(K)+4\Arf(K)$ where it is fixed, also excludes exactly the knots that the linking form and
Proposition~\ref{prop:dtest} exclude (\S\ref{sec:controls}).  So Proposition~\ref{prop:donaldson} adds nothing for these alternating knots; but it needs no correction terms, so it applies to non-alternating knots whose correction terms we
cannot compute, as soon as some diagram has a checkerboard surface that is definite of the right sign
(\S\ref{sec:lower}).

\subsection{Surfaces built from band moves}\label{sec:bandconv}
\begin{figure}[ht]
\centering
\includegraphics{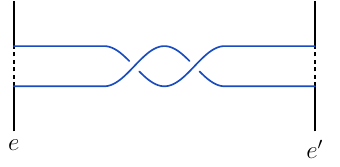}
\caption{A band joining the edges $e$ and $e'$ of a face of a knot diagram, with $t=2$: its long sides cross twice, and
at each crossing the strand from lower left to upper right is over, that is, the long side on the right passes over the
one on the left as one walks along the band, in either direction; for $t<0$ the other strand is over at each of the
$|t|$ crossings.  The band move replaces the two short arcs of $e$ and
$e'$ between the long sides, drawn dashed, by the long sides.}
\label{fig:band}
\end{figure}

We describe band moves on knot diagrams as in Figure~\ref{fig:band}.  A band move on a knot gives a knot or a
two-component link; a band move between two knots is nonorientable, since a coherent band would give a link of two
components.

If a band move turns $K$ into a knot $J$ and $J$ bounds a surface $F'\subset B^4$, then the band, pushed into a collar,
together with $F'$ is a surface $F$ with boundary $K$; it is nonorientable, and $b_1(F)=b_1(F')+1$.  So a slice disc of
$J$ gives a M\"obius band, a M\"obius band of $J$ a punctured Klein bottle, and a locally flat disc a locally flat
M\"obius band.  We obtain the surfaces of $J$ from further band moves.  If $J$ is the unknot, $F'$ is a disc; a band from
$J$ to the two-component unlink, or two bands to the three-component unlink, give a ribbon disc made of two discs and
the band, or of three discs and the two bands; and a band from $J$ to the unknot or to such a ribbon knot gives a
M\"obius band of $J$.  The unknot and the unlinks are recognised by Reidemeister moves to diagrams
without crossings.  So each smooth surface that we build is given by a sequence of band moves on PD codes, some of them related by
Reidemeister moves, that ends in a diagram which Reidemeister moves reduce to an unlink, and it depends neither on an
identification of knots nor on any value of \KI{}.

\section{Computations}\label{sec:computations}

We applied the bounds of Section~\ref{sec:tools} to the prime knots with at most $13$ crossings, and the linking-form
bound and Proposition~\ref{prop:dtest} to those with $14$ crossings (\S\ref{sec:fourteen}).  All homology groups,
linking forms, Jordan decompositions and lattices are computed exactly.  The lower bounds are in \S\ref{sec:lower}, the
surfaces that give the upper bounds in \S\ref{sec:upper}, the resulting values in \S\ref{sec:values} and the consistency checks in \S\ref{sec:controls}.

\subsection{The knot tables}\label{sec:methods}

For the $\nKnotsThirteen$ prime knots with at most $13$ crossings we use \KI{} \cite{KnotInfo}, version 2026.11.2: the
values and ranges of $\gfour$, and the Seifert matrices, PD codes and Montesinos notation; versions 2026.8.1, 2026.10.5
and 2026.11.2 of \KI{} agree on all invariants and diagrams that we use.  \KI{} gives the least first Betti number of a nonorientable surface in $B^4$,
which is $1$ for slice knots; for knots that are not slice it is $\gfour$.  It lists a range for $\nOpen$ of these
knots: $[1,2]$ for $\nOpenOneTwo$, $[1,3]$ for $\nOpenOneThree$, $[2,3]$ for $\nOpenTwoThree$, $[2,4]$ for
$\nOpenTwoFour$, $[1,4]$ for $\nOpenOneFour$ knots and $[1,5]$ for one knot.  Most of its values for $11$ to $13$
crossings come from the work of Lee and Sabloff \cite{LS}, who find M\"obius bands and obstruct them with
Yasuhara's condition and with the criterion of \cite[Cor.~2.7]{MY}, stated as \cite[Thm.~2.4]{LS}, when the determinant
is square-free.  For the $\nLS$ knots with $14$ crossings we use the table of \cite{LS} and the diagrams of SnapPy
\cite{SnapPy}.

\subsection{Lower bounds}\label{sec:lower}
Proposition~\ref{prop:lower} collects the lower bounds from the linking form and from the correction terms of Goeritz
matrices and star-shaped plumbings, and Proposition~\ref{prop:newd} those from plumbing trees, surgery descriptions and
Donaldson's theorem for the non-alternating knots that are not Montesinos.

\begin{proposition}\label{prop:lower}
Among the $\nKnotsThirteen$ prime knots with at most $13$ crossings:
\begin{enumerate}[label=(\alph*)]
\item $h_\lambda(K)\le3$ and $h^\oplus(K)=h_\lambda(K)$ for all of them, and $h_\lambda(K)=3$ for $\nHthree$ of them,
which therefore have $\gtop(K)\ge3$;
\item of the $\nRankOneTargets$ knots whose \KI{} range starts at $1$, $\nRankOneObstr$ have $h_\lambda\ge2$
($\nRankOneObstrOneTwo$ with range $[1,2]$ and $\nRankOneObstrOneThree$ with range $[1,3]$), and so does every knot with
range $[2,3]$ or $[2,4]$; all of these have $\gtop(K)\ge2$;
\item of the other knots whose range starts at $1$, Proposition~\ref{prop:dtest} shows $\gfour(K)\ge2$ for
$\nDAltObstr$ of the $\nDAltTargets$ alternating knots ($\nDAltObstrOneTwo$ with range $[1,2]$ and
$\nDAltObstrOneThree$ with range $[1,3]$) and for $\nDMontObstr$ of the $\nDMontTargets$ non-alternating Montesinos
knots ($\nDMontObstrOneTwo$ with range $[1,2]$ and $\nDMontObstrOneThree$ with range $[1,3]$).
\end{enumerate}
\end{proposition}

\begin{proof}
(a) is Theorem~\ref{thm:compute} and Proposition~\ref{prop:hoplus}, since every determinant is at most
$\nMaxDetThirteen<3^9$, together with Theorem~\ref{thm:lower}; (b) is Corollary~\ref{cor:criteria}(b) and
Theorem~\ref{thm:lower}; (c) is Proposition~\ref{prop:dtest}, with the correction terms of the Goeritz manifold or of the
star-shaped plumbing (\S\ref{sec:dtest}).
\end{proof}

\para{Computation.}
The form $\lambda_K$ is computed from $A=V+V^{T}$: a Smith normal form gives $H_K$, and for each prime $p\mid D$ a
diagonalisation of $A$ over $\Z_{(p)}$ gives the Jordan decomposition of $\lambda_{K,p}$
\cite[(3)]{LeeTheory}; Theorem~\ref{thm:compute} then gives $h_\lambda(K)$.  For an alternating knot we use the Goeritz matrix of the \KI{} diagram, which is alternating, and for a Montesinos
knot the star-shaped plumbing given by its Montesinos notation, as in \cite{Lee}; for the other non-alternating
knots we use the methods of Proposition~\ref{prop:newd}.  For Goeritz matrices and plumbings the
correction terms are maxima over characteristic covectors (\S\ref{sec:dtest}), computed exactly.  The test of Proposition~\ref{prop:dtest} covers both orientations $Z=\pm\Sigma_K$, all $\ell$ with $D/\ell$ a square, all
isotropic $M$ with $\ell|M|^2=D$ and all admissible $g$.

\begin{proposition}\label{prop:newd}
The following $\nNewDAll$ non-alternating knots $K$ with at most $13$ crossings that are not Montesinos, have
$h_\lambda(K)=1$ and have the \KI{} range $[1,2]$ or $[1,3]$ satisfy $\gfour(K)\ge2$:
\begin{enumerate}[label=(\alph*)]
\item $\nNDPlumbObstr$ knots for which $\Sigma_K$ or $-\Sigma_K$ bounds a negative-definite plumbing tree with at most
two bad vertices whose correction terms exclude a smooth M\"obius band;
\item $\nSurgExcl$ knots with unknotting number one for which a surgery description of $\Sigma_K$ excludes a smooth
M\"obius band, $\nSurgExclOnly$ of them not among the knots of (a);
\item $\nDefCert$ further knots with a diagram that has a checkerboard surface $V$ to which
Proposition~\ref{prop:donaldson} applies and for which the lattice $(H_1(V),-\varepsilon\mathcal G_V)$, of rank $8$ or
$9$, does not embed in $\Z^{b_1(V)+1}$.
\end{enumerate}
\end{proposition}

\begin{proof}
In (a) and (b) Proposition~\ref{prop:dtest} applies, with the correction terms of the plumbing tree
\cite[Cor.~1.5]{OSz-plumb} or of the surgery description, bounded or computed as in \S\ref{sec:dtest}; in (c)
Proposition~\ref{prop:donaldson} applies.
\end{proof}

We applied (a) and (b) to the $\nNDTargets$ non-alternating knots with at most $13$ crossings that are not Montesinos,
have $h_\lambda\le1$ and the \KI{} range $[1,2]$ or $[1,3]$, and bound none of the M\"obius bands of the first
construction of \S\ref{sec:upper}, and (c) to the $\nDefKnots$ of them that (a), (b) and the second construction of
\S\ref{sec:upper} leave open.  The knots of the proposition are marked in Appendices~\ref{app:newtwo}
and~\ref{app:newthree}.

\para{Plumbing trees.}
For each of the $\nNDTargets$ knots we decided whether $\Sigma_K$ or $-\Sigma_K$ bounds a negative-definite plumbing
tree with the method of \cite{LeeNDP}, whose results already decide this for $\nNDFollowup$ of them.  The double branched cover is built in SnapPy as the double cover of the orbifold obtained by
$(2,0)$-filling the exterior of $K$, and recognised with Regina \cite{Regina}.  When it is a graph manifold, the Neumann
normal forms of $\Sigma_K$ and $-\Sigma_K$, computed from the recognised Seifert pieces or from an algebraic (arborescent) diagram of
$K$ \cite[Lemma~2.7]{LeeNDP}, decide the question: every normal form that we computed is a tree, and a graph manifold that is a rational homology sphere and whose normal
form is a tree bounds a negative-definite plumbing tree if and only if its normal form is negative definite
\cite[Cor.~2.5]{LeeNDP}.  Each tree that we use was identified with $\pm\Sigma_K$ by triangulating its plumbed manifold,
and has $|\det Q|=D$.  Of the $\nNDTargets$ covers, $\nNDHyp$ are hyperbolic, as verified with interval arithmetic,
$\nNDMixed$ are toroidal with a hyperbolic JSJ piece, $\nNDNoDefAll$ are graph manifolds of which neither orientation
bounds a negative-definite plumbing tree, and $\nNDNegDef$ bound one; $\nNDBadTwo$ of their trees have two bad
vertices and the others one.  The test of Proposition~\ref{prop:dtest} passes only for $\kn{13n}{2413}$ and
$\kn{13n}{3505}$.  An independent calculation for cyclic $H_K$ agrees on the $\nNDIndependent$ knots to which
it applies, and for $\nNDStarRoute$ of the knots a star-shaped plumbing built directly from the Seifert invariants that
Regina recognises gives the same correction terms.

\para{Surgery descriptions.}
For each of the $\nSurgTargets$ of the $\nNDTargets$ knots that have unknotting number one we found a crossing of the
\KI{} diagram, or of a diagram related to it by
Reidemeister moves, whose change gives a diagram that Reidemeister moves reduce to the unknot $U$.  A crossing circle
$c$ at this crossing bounds a disc that meets $U$ in two points, and a $(\pm1)$-surgery on $c$ turns $U$ into $K$.  The
two lifts of $c$ to $\Sigma_U=S^3$ are parallel copies of the lift $C$ of the crossing arc, and $\Sigma_K$ is a surgery
on $C$ with slope $\pm D/2$, as in the Montesinos trick.  With SnapPy we built the double cover of the orbifold obtained
by $(2,0)$-filling $U$ in the exterior of $U\cup c$, found with Regina the filling of the two lifts that gives $S^3$,
and so obtained the exterior of $C$ with its meridian.  We then found the slope $eD/2$, $e=\pm1$, for which
$S^3_{eD/2}(C)$ is homeomorphic to $\Sigma_K$: by an isometry of hyperbolic structures, by drilling the core of the
filling from $\Sigma_K$, or by Regina's recognition.  The sign $e$ was confirmed by the linking form and by the
Casson--Walker invariant, with Mullins' formula as in \cite[\S3.4]{LeeNDP}.  A diagram of $C$ comes from SnapPy, and its
chirality is fixed by an orientation-preserving isometry of exteriors, or by Regina's names of the surgeries on the
diagram and on its mirror image.

The two bounds of \S\ref{sec:dtest} exclude a smooth M\"obius band for
$\nSurgBox$ of these $\nSurgTargets$ knots.  For the $\nSurgExactTargets$ of them that neither these bounds nor (a) exclude, we computed
the exact
$V_j(C)$ as in
\S\ref{sec:dtest}; for each of them both ways determine the $V_s(C)$, and they agree.  The chirality of $C$ was fixed as above, or by the reduced Khovanov homology of $K$: if its
rank over $\Z/2$, as listed by \KI{}, is $D$, then $\Sigma_K$ is an $L$-space \cite[Cor.~1.2]{OSz-bdc}, and only one
chirality of $C$ has an $L$-space surgery of slope $D/2$.  For the other knots we compared the connected covers of
small degree, with their first homology, of $\Sigma_K$ and of the two possible surgeries, whose covers differ, with SnapPy \cite{SnapPy}.  With the exact correction terms,
Proposition~\ref{prop:dtest} excludes a smooth M\"obius band for $\kn{13n}{3844}$ and $\kn{13n}{3855}$ and passes
for the other $\nSurgExactPass$ knots.  For these two knots the covers fixed the chirality of $C$, and once it is fixed,
the bounds above exclude them as well.

\para{The other knots.}
Of the $\nNDTargets$ knots, $\nNDLeft$ remain open after (a) and (b).  For $\nNDLeftNoMethod$ of them neither $\Sigma_K$
nor $-\Sigma_K$ bounds a negative-definite plumbing tree and \KI{} does not list unknotting number one; for
$\kn{13n}{2413}$ and $\kn{13n}{3505}$ the test passes; and for the other $\nSurgLeft$ knots with unknotting number one it
passes with the exact correction terms.  The method of (b) also applies to $\kn{11n}{176}$ (\S\ref{sec:corrections}).
The second construction of \S\ref{sec:upper} gives M\"obius bands for four of these knots, which leaves the $\nDefKnots$ knots
of (c).

\para{Definite surfaces.}
For (c), the congruence $\sigma(K)+4\Arf(K)\equiv\pm2\pmod 8$ holds for $\nDefSearched$ of the $\nDefKnots$ knots; for the other $\nDefSZero$ it is $0$ modulo $8$, and Proposition~\ref{prop:donaldson} does not apply.  For none of these knots does the \KI{} diagram have a checkerboard surface with a definite form, whereas non-minimal diagrams often do.  For each of the $\nDefSearched$ knots
we examined the \KI{} diagram, the diagrams of the two constructions of M\"obius bands in \S\ref{sec:upper}, and
diagrams made from the \KI{} diagram by random sequences of Reidemeister moves, with at most ten crossings more;
$\nDefDiagrams$ distinct diagrams in all.  For each checkerboard surface whose form is definite of the required sign,
checking all vectors of the required norm decides whether the lattice embeds in the standard lattice of rank one more.  Some diagram has
such a surface for $\nDefWithSurface$ of the knots, and for $\nDefCert$ of them one of these lattices does not embed; the diagrams that carry these lattices have $\nDefCertMinC$ to $\nDefCertMaxC$ crossings.

\para{Verification.}
Each of these obstructions was checked independently.  An orientation-preserving isometry from the exterior of the diagram to that of the \KI{} diagram takes the meridian to the meridian, confirming that the diagram represents $K$ rather than its mirror image.  The Goeritz matrices recomputed from the diagram give the same positive-definite lattice of determinant $\det K$ up to isometry and the sign $-\varepsilon$, and this lattice presents the linking form of
$-\varepsilon$ times the symmetrised Seifert matrix of \KI{}; this compares the Goeritz matrix with the Seifert matrix,
including the sign, since for each of the $\nDefCert$ knots this linking form is not isometric to its negative.  Finally, a separate check confirms that the lattice does not embed.

\para{Surgery descriptions by drilling.}
For the $\nDrillTargets$ of the $\nDefKnots$ knots to which (c) was applied whose double branched cover is hyperbolic and whose unknotting number is not one in \KI{}, we also looked for surgery
descriptions without an unknotting crossing.  We drilled $\Sigma_K$ along short closed curves; when the drilled manifold has a filling that Regina recognises as $S^3$, it is
the exterior of a knot $C\subset S^3$, and $\Sigma_K\cong S^3_{p/q}(C)$.  The slope, the chirality of $C$ and its
correction terms were then treated as in (b).  This gave surgery descriptions
for $\nDrillDescr$ knots.  For $\nDrillDetermined$ of them the test of Proposition~\ref{prop:dtest} is decided, with
$V_j(C)$ read off the knot Floer homology of $C$ when $C$ or its mirror image is an $L$-space knot or $C$ is thin, and
otherwise from the knot Floer complex of $C$ as in (b); for the other $\nDrillNoHFK$ we could not compute the knot Floer homology of $C$.
Proposition~\ref{prop:dtest} excludes a smooth M\"obius band for $\kn{13n}{2626}$, $\kn{13n}{2632}$ and
$\kn{13n}{3616}$, each with $q=3$ and $C$ an $L$-space knot, and passes for the other $\nDrillPass$.  The three knots are among the $\nDefCert$ knots of (c), so for them two independent methods give $\gfour\ge2$.

\subsection{Upper bounds}\label{sec:upper}
The upper bounds come from surfaces built from band moves as in \S\ref{sec:bandconv}, in two constructions that are
described after Proposition~\ref{prop:upper}.

\begin{proposition}\label{prop:upper}
Among the prime knots with at most $13$ crossings for which \KI{} lists a range of $\gfour$:
\begin{enumerate}[label=(\alph*)]
\item $\nGammaOneAll$ bound smooth M\"obius bands built from band moves: the four knots $\kn{13a}{708}$,
$\kn{13a}{3326}$, $\kn{13a}{3609}$ and $\kn{13a}{4877}$, whose range is $[1,4]$ or $[1,5]$, $\nMobiusOneTwoAll$ whose
range is $[1,2]$ and $\nMobiusOneThreeAll$ whose range is $[1,3]$;
\item $\nKleinTotal$ bound smooth punctured Klein bottles built from band moves: every knot of range $[1,3]$, $[2,3]$ or
$[2,4]$ with $h_\lambda\le2$ that Proposition~\ref{prop:lower} obstructs, that is, $\nKleinCertified$ of range $[1,3]$
and the $\nRangesTargets$ of range $[2,3]$ or $[2,4]$, which have $h_\lambda=2$, and $\nMobiusKlein$ of the knots of range
$[1,3]$ that it does not obstruct;
\item $\kn{13a}{2890}$, whose range is $[2,4]$, bounds a smooth nonorientable surface with $b_1=3$ built from band moves.
\end{enumerate}
\end{proposition}

\begin{proof}
The surfaces are built as in \S\ref{sec:bandconv}, from the bands described below.
\end{proof}

\para{The first construction.}
We try the bands in the faces of the \KI{} diagram and of further diagrams related to it by Reidemeister
moves.  Each knot of (a) found in this way bounds a M\"obius band made of one band and a
disc: a band turns it into a ribbon knot $J'$, whose ribbon disc it bounds with the band, or, for $\nMobiusUnknot$ knot,
into the unknot.  The band lies in the \KI{} diagram for three of the knots of range $[1,4]$ or $[1,5]$ and for
$\nMobiusKI$ of the $\nMobiusAll$ knots of range $[1,2]$ or $[1,3]$, and in a diagram related to it by Reidemeister moves
for $\kn{13a}{4877}$ and the other $\nMobiusStageB$.  In the \KI{} diagram of each of the $\nKleinCertified$ knots of (b) with range $[1,3]$ that
Proposition~\ref{prop:lower} obstructs, we found a band of Figure~\ref{fig:band}, with $t=0$ for $\nKleinTZero$ of them
and $|t|=1$ for the others, that leads to a punctured Klein bottle:
$J$ has a band to the unknot for $\nKleinUnknotBand$ of them and a band to a ribbon knot for the other
$\nKleinRibbonBand$.  The same construction gives such a surface for each of the $\nRangesTargets$ knots of (b) with range
$[2,3]$ or $[2,4]$, with $t=0$ for $\nRangesTZero$ of them and $|t|=1$ for the others; $J$ has a band to the unknot for
$\nRangesUnknotBand$ of them and a band to a ribbon knot for the other $\nRangesRibbonBand$.  For $\nRangesKI$ of these
knots the first band lies in the \KI{} diagram, and for the other $\nRangesStageB$ in a diagram related to it by
Reidemeister moves; for one of the latter, $\kn{13a}{3126}$, the bands of $J$ lie in a diagram obtained by Reidemeister
moves from the diagram of $J$ that the first band produces.  The last $\nMobiusKlein$ knots of (b) are the knots of range $[1,3]$ that
pass the tests of Proposition~\ref{prop:lower} and for which this construction gives no M\"obius band; their punctured Klein
bottles lie in the \KI{} diagram for $\nMobiusKleinKI$ of them.  For $\kn{13a}{2890}$ a band with $t=1$ gives
$\kn{7}{5}$, and two further bands turn $\kn{7}{5}$ into the unknot; with a disc, the three bands form a nonorientable
surface with $b_1=3$.  Every surface used in the proofs was reconstructed independently from the \KI{} PD codes: the final unknot or unlink is recognised by Reidemeister moves alone, and the further
diagrams of $K$, as well as the diagram of $J$ used for $\kn{13a}{3126}$, are related by Reidemeister
moves to the diagrams from which they come.

\para{The second construction.}
For the $\nWideTargets$ knots of range $[1,2]$ or $[1,3]$ that the first construction and Propositions~\ref{prop:lower} and~\ref{prop:newd}(a), (b) leave open, we looked again for M\"obius bands made of band moves.  When SnapPy identifies the knot $J$ produced by the first band with a knot that \KI{} lists as
slice, we look for a ribbon disc of $J$, made of one band to the two-component unlink or of two bands to the
three-component unlink, also in diagrams of $J$ obtained by Reidemeister moves from the diagram that the first band
produces, and in diagrams of the identified knot or of its mirror image that Reidemeister moves relate to that diagram;
the identification serves only to choose where to look.  This gives M\"obius bands for $\nWideMobius$ knots, $\nWideMobiusAlt$ of them
alternating.  For all $\nWideMobius$ the first band lies in a diagram that the first construction had already examined, the \KI{} diagram for $\nWideKI$ of them, and the ribbon disc lies
in another diagram of $J$: one obtained by Reidemeister moves for $\nWideReached$ of them, and a diagram of the
identified knot or of its mirror image for the other $\nWideViaX$.  The ribbon disc is made of one band and two discs,
except for $\kn{13n}{1374}$, where it is made of two bands and three discs.  For the knots still left open after Proposition~\ref{prop:newd}(c), we searched in the same way for a ribbon disc for every knot $J$ produced by these bands.  For each knot $K$, the resulting knots
$J$ were grouped by their determinant, Alexander polynomial and hyperbolic volume, and the groups were
ordered and filtered by necessary conditions for sliceness computed from a Seifert matrix and from knot Floer homology,
which again serve only to choose where to look; for each group we looked for a ribbon disc in the diagram of $J$ that
the first band produces and, made of one band, in up to $300$ diagrams of $J$ obtained from it by Reidemeister moves.
This gives M\"obius bands for four more knots, $\kn{13a}{1412}$, $\kn{13a}{2898}$, $\kn{13a}{3125}$ and
$\kn{13n}{3901}$, each with a ribbon disc made of one band and two discs; for $\kn{13a}{1412}$ the knot $J$ is the one
that SnapPy identifies with $\kn{12a}{77}$, for which the search above had found no ribbon disc in a diagram related to
that of $J$.  For $\nJStopped$ of the knots left open, not all groups were
examined.  As in the first construction, every band was reconstructed from the \KI{} PD codes, together with Reidemeister moves from the \KI{} diagram
of $K$ to the diagram of the first band and from the resulting diagram of $J$ to that of the
ribbon disc; thus these M\"obius bands depend neither on the identification of $J$ nor on any value in \KI{}.  Besides the diagrams of the first construction we used new ones, made by longer sequences of random
Reidemeister moves; they have at most ten crossings more than the crossing number of $K$, and SnapPy verified that each
has the exterior of the \KI{} knot.  Each knot for which no M\"obius band was found was examined in
$\nWideOpenMin$ to $\nWideOpenMax$ diagrams; $\nWideOpenNew$ of these diagrams were new.

\subsection{The values}\label{sec:values}

Together with the bounds that \KI{} lists, Propositions~\ref{prop:lower}--\ref{prop:upper} determine $\gfour$ for
$\nDetermined$ of the $\nOpen$ knots with a range, as follows; Table~\ref{tab:ranges} summarises the effect on the
ranges.

\para{$\gfour=3$.}
The $\nGammaThreeAll$ knots with $h_\lambda=3$ and range $[2,3]$ ($\nThreeFromTwoThree$ knots), $[1,3]$
($\nThreeFromOneThree$ knots) or $[2,4]$ (the knot $\kn{13a}{2890}$) have $\gfour=\gtop=3$, by
Proposition~\ref{prop:lower}(a) and the upper bound $3$ of \KI{} or of Proposition~\ref{prop:upper}(c).

\para{$\gfour=2$.}
The knots that Propositions~\ref{prop:lower}(b), (c) and~\ref{prop:newd} obstruct and that have $h_\lambda\le2$ have
$\gfour=2$, with the upper bound $2$ of their range $[1,2]$ or of the punctured Klein bottles of
Proposition~\ref{prop:upper}(b): $\nGammaTwoOneTwo$ knots of range $[1,2]$, $\nGammaTwoOneThree$ of range $[1,3]$,
$\nRangesTwoThree$ of range $[2,3]$ and $\nRangesTwoFour$ of range $[2,4]$, $\nGammaTwoAll$ in all.  The $\nGammaTopTwo$
of them with $h_\lambda=2$ have $\gtop=2$ as well.

\para{$\gfour=1$.}
The $\nGammaOneAll$ knots of Proposition~\ref{prop:upper}(a) have $\gfour=1$, since none of them is smoothly slice.  The
four knots of range $[1,4]$ or $[1,5]$ have signatures $-6$, $-8$, $-4$ and $-4$.  Of the $\nMobiusAll$ knots of range
$[1,2]$ or $[1,3]$ found by the first construction, $\nMobiusHOne$ have $h_\lambda=1$, $\nMobiusSig$ have a nonzero signature,
and the remaining $\nMobiusKIgfour$ have smooth $4$-genus $1$ in \KI{}; of the $\nWideAll$ found by the second construction, $\nWideHOne$ have $h_\lambda=1$, $\kn{13a}{3125}$ and $\kn{13a}{4799}$ have nonzero signatures, and the other $\nWideKIgfour$ have
smooth $4$-genus $1$ in \KI{}.  Those with $h_\lambda=1$ or a nonzero signature are not topologically slice either.
Neither are $\nMobiusFoxMilnor$ of the $\nMobiusKIgfour$, the knots $\kn{12n}{837}$, $\kn{13a}{1431}$, $\kn{13a}{1506}$,
$\kn{13n}{1967}$ and $\kn{13n}{2985}$, nor $\kn{13a}{4729}$ among the $\nWideKIgfour$.  By the theorem of Fox and
Milnor (\cite[Thm.~2]{FoxMilnor} for piecewise linear discs, \cite[Thm.~15.2]{FNOP} for locally flat ones), the
Alexander polynomial of a topologically slice knot has the form
\[
\Delta(t)=\pm t^nf(t)f(t^{-1}),\qquad f\in\Z[t],
\]
so each of its irreducible factors $g$ with $t^{\deg g}g(t^{-1})=\pm g(t)$ has even multiplicity; but every irreducible
factor of the Alexander polynomials of these six knots has this symmetry and multiplicity~$1$.  So $\gtop=1$ for
$\nGammaTopOneAll$ knots in all.  The Alexander polynomials of the other $\nMobiusFoxMilnorHolds$ and $\nWideFMHolds$ have the form
$\pm t^nf(t)f(t^{-1})$.

\para{The ranges.}
So every range with lower bound $2$ is settled, and every range $[1,3]$ that is not settled becomes $[1,2]$; the range $[1,2]$ remains
for $\nOpenLeft$ knots, $\nOpenLeftFromOneThree$ of which had the range $[1,3]$.  Of these $\nOpenLeft$ knots,
$\nOpenLeftHOne$ have $h_\lambda=1$, and the remaining one, $\kn{13n}{3505}$, has smooth $4$-genus $1$ in \KI{}; so none
of them is smoothly slice, and their range is $[1,2]$.
They are listed in Appendix~\ref{app:openleft}, and the knots with $\gfour=2$ and range $[1,2]$ or $[1,3]$ in
Appendices~\ref{app:newtwo} and~\ref{app:newthree}; every other knot of range $[1,2]$ or $[1,3]$ has $\gfour=1$, except the
$\nThreeFromOneThree$ knots of range $[1,3]$ with $\gfour=3$.  The knots with $\gfour=2$ and
range $[2,3]$ or $[2,4]$ are all those with these ranges and $h_\lambda=2$.

\para{Sharpness of $h_\lambda$.}
The $\nHthree$ knots with $h_\lambda=3$ are the $\nGammaThreeAll$ above, $\kn{8}{18}$, $\kn{9}{49}$ and $\kn{11n}{133}$;
for $\kn{8}{18}$ \KI{} lists $3$, and the other two are discussed in \S\ref{sec:corrections}.  Table~\ref{tab:hthree}
lists them with the invariants that give $h_\lambda=3$.  So each of the $\nSharpAll$ prime knots with at most $13$
crossings and $h_\lambda\ge2$ has $\gfour=\gtop=h_\lambda$: those with $h_\lambda=3$ by the above, by
\S\ref{sec:corrections} and by the value of $\kn{8}{18}$, and those with $h_\lambda=2$ by the above or by the value $2$
that \KI{} lists for the others.  The lower bound $h_\lambda$ is therefore sharp on all of them.

\begin{table}[ht]
\centering
\footnotesize\setlength{\tabcolsep}{1.5pt}
\begin{tabular}{lrrrrrrrrrrr}
\toprule
\KI{} range & knots & $h_\lambda\ge2$ & \ref{prop:lower}(c) & \ref{prop:newd}(a,b) & \ref{prop:newd}(c) &
$h_\lambda=3$ & bands & $\gfour=1$ & $\gfour=2$ & $\gfour=3$ & open \\
\midrule
$[1,2]$ & $\nOpenOneTwo$ & $\nRankOneObstrOneTwo$ & $\nDObstrOneTwo$ & $\nNDOneTwo$ & $\nDefOneTwo$ & &
$\nMobiusOneTwoAll$ & $\nMobiusOneTwoAll$ & $\nGammaTwoOneTwo$ & & $\nOpenLeftOneTwo$ \\
$[1,3]$ & $\nOpenOneThree$ & $\nRankOneObstrOneThree$ & $\nDObstrOneThree$ & $\nNDOneThree$ & $\nDefOneThree$ &
$\nThreeFromOneThree$ & $\nBandsOneThree$ & $\nMobiusOneThreeAll$ & $\nGammaTwoOneThree$ & $\nThreeFromOneThree$ &
$\nOpenLeftFromOneThree$ \\
$[2,3]$ & $\nOpenTwoThree$ & & & & & $\nThreeFromTwoThree$ & $\nRangesTwoThree$ & & $\nRangesTwoThree$ &
$\nThreeFromTwoThree$ & \\
$[2,4]$ & $\nOpenTwoFour$ & & & & & $1$ & $\nOpenTwoFour$ & & $\nRangesTwoFour$ & $1$ & \\
$[1,4]$, $[1,5]$ & $\nOpenOneFourFive$ & $0$ & $0$ & $0$ & $0$ & & $\nOpenOneFourFive$ & $\nOpenOneFourFive$ & & & \\
\midrule
all & $\nOpen$ & & & & & & & $\nGammaOneAll$ & $\nGammaTwoAll$ & $\nGammaThreeAll$ & $\nOpenLeft$ \\
\bottomrule
\end{tabular}
\caption{The effect of Propositions~\ref{prop:lower}--\ref{prop:upper} on the \KI{} ranges.  The columns $h_\lambda\ge2$,
\ref{prop:lower}(c), \ref{prop:newd}(a,b) and \ref{prop:newd}(c) count the knots excluded by the linking form
(Proposition~\ref{prop:lower}(b)), and, among those that pass it, the knots excluded by the correction terms of Goeritz
matrices and star-shaped plumbings (Proposition~\ref{prop:lower}(c)), by those of plumbing trees and surgery descriptions
(Proposition~\ref{prop:newd}(a), (b)) and by Donaldson's theorem (Proposition~\ref{prop:newd}(c)).  These sets are
disjoint, and all of their knots have $\gfour\ge2$; for the range $[1,3]$ the first two columns count
$\nRankOneObstrOneThree+\nDObstrOneThree=\nNewTwoThree$ knots, and the $\nThreeFromOneThree$ with $h_\lambda=3$ among them
have $\gfour=3$.  For the ranges $[2,3]$ and $[2,4]$ these columns are left empty: the lower bound is already $2$, and
every such knot has $h_\lambda\ge2$.  The column bands counts the knots whose upper bound comes from the surfaces of
Proposition~\ref{prop:upper}: M\"obius bands for the knots counted under $\gfour=1$, and punctured Klein bottles or a
surface with $b_1=3$ for the others, among them $\nKleinFromOneThreeLeft$ of the $\nFirstOpenFromOneThree$ knots of range $[1,3]$ that
Proposition~\ref{prop:lower} does not obstruct and for which the first construction gives no M\"obius band (the other
$\nWideFromOneThree$ bound M\"obius bands of the second construction);
Proposition~\ref{prop:newd} settles $\nNDOneThree$ of them by (a), (b) and $\nDefOneThree$ by (c).  The columns $\gfour=1$, $2$ and $3$ count the values of
\S\ref{sec:values}, $\nDetermined$ in all, and the column open counts the knots whose value remains open; their range is
now $[1,2]$.}
\label{tab:ranges}
\end{table}

\subsection{Consistency checks}\label{sec:controls}
The following computations compare the methods of \S\S\ref{sec:lower}--\ref{sec:upper} with one another and with the
knots whose value is known.

\para{Linking forms.}
We verified that both congruences of \cite[Lemma~2.5]{LeeTheory} hold for all $\nKnotsThirteen$ knots.  For the
$\nGoeritzChecked$ knots to which these tests are applied, here or in \S\ref{sec:lower}, $\lambda_K$ was recomputed from both Goeritz matrices
of the \KI{} diagram, with the same result up to one global sign.  Two further computations verify the cases $h_\lambda\ge2$ and
$h_\lambda\ge3$.  The first decides whether $\lambda_K\cong\langle\tau/\ell\rangle\oplus\mu$ with $\mu$ metabolic for
some $\ell$ and $\tau$; by Proposition~\ref{prop:hoplus} this is equivalent to the Witt-level condition when no prime is
exceptional. It agrees with a direct check of all elements of $H_K$ for the $\nValidate$ knots of determinant less than
$200$.  The second examines all classes of odd binary forms $Q$ of determinant $\pm\ell$ with $D/\ell$ a square and
decides, prime by prime, whether $\lambda_Q$ is an orthogonal summand of $\lambda_K$ with metabolic complement.

\para{Known values.}
The tests exclude no knot that \KI{} 2026.11.2 lists with the value $1$; for those of Proposition~\ref{prop:newd} see
below.  The only knot that they exclude although \cite{Fairchild}, the table of \cite{LS} or an earlier
version of \KI{} gives $\gfour=1$ is $\kn{11n}{57}$ (\S\ref{sec:corrections}).  Among the $\nControlsOne$ knots with at
most $13$ crossings for which \KI{} lists the value $1$, none has $h_\lambda\ge2$; Proposition~\ref{prop:dtest} excludes
none of the $\nDAltControls$ alternating ones and none of the $\nDMontControls$ non-alternating Montesinos ones.  Among the
$\nRankTwoControls$ knots for which \KI{} lists an upper bound at most $2$, exactly two have $h_\lambda=3$, namely
$\kn{9}{49}$ and $\kn{11n}{133}$ (\S\ref{sec:corrections}).  Conversely, the first construction gives no M\"obius band for the $\nMobiusControls$ knots of range $[1,2]$ or
$[1,3]$ that one of the tests excludes, in $\nMobiusControlDiagrams$ diagrams, and the $\nHthree$ knots with
$h_\lambda=3$ showed no band to a knot with \KI{} value $1$ in $\nKleinControlDiagrams$ diagrams.  The second construction gives no M\"obius band and no band to a knot identified with a slice knot of \KI{} for $\nWideControls$ knots for which
Propositions~\ref{prop:lower} and~\ref{prop:newd} give $\gfour\ge2$, in $\nWideControlDiagrams$ diagrams, and its
continuation to every knot $J$ gives no M\"obius band for $\nJControls$ of them.

\para{Plumbing trees and surgeries.}
The test with plumbing trees excludes no knot with $\gfour=1$ (in \KI{} or by Proposition~\ref{prop:upper}(a)): it passes for all $\nNDControlsOne$ such knots among the $\nNDControlsFollowup$ knots with a
negative-definite tree in the tables of \cite{LeeNDP}, as well as for all $\nNDControlsTrees$ negative-definite trees that the method finds for $\nNDControlsTreeKnots$ of the $\nNDControlsPipeline$ non-alternating, non-Montesinos knots with $\gfour=1$ outside the tables of \cite{LeeNDP}.
Among the knots tested in Proposition~\ref{prop:lower}(c) and in these checks, $\nNDControlsKnown$ alternating and Montesinos
knots have an algebraic \KI{} diagram, in the sense of \cite[Lemma~2.7]{LeeNDP}; for $\nNDControlsKnownTested$ of them the method finds negative-definite trees, and their correction terms agree with those of the Goeritz matrices and
star-shaped plumbings, with the same conclusions.  The surgery method was applied to $\nSurgCalib$ of the alternating knots
with \KI{} unknotting number one that are not Montesinos and have hyperbolic double branched covers.  For $\nSurgCalibAgree$
of them the correction terms of the surgery description, with $V(C)$ read off the knot Floer complex of $C$, agree with
those of the Goeritz manifold up to an isometry of linking forms; for the remaining one, $\kn{13a}{1085}$, we could not compute the knot Floer complex of $C$.  The two bounds alone exclude $\nSurgCalibBox$ of the $\nSurgCalib$ knots, exactly those that these correction
terms exclude, and none with \KI{} value $1$.  The computation of the exact $V_j(C)$ reproduces these correction terms and conclusions, excludes again the $\nSurgExactControls$ knots with unknotting number one that Proposition~\ref{prop:newd}(a) and the
bounds of \S\ref{sec:dtest} exclude,
gives the $V_j(C)$ of $\kn{11n}{176}$ in \S\ref{sec:corrections}, and gives the $V_s$ of the tensor product of the staircase
complexes of the two summands for $\nSurgTensor$ connected sums of two torus knots.  The labelling of the Spin$^c$ structures in the formula of Ni and Wu
(\S\ref{sec:dtest}) was verified on lens spaces against their linear plumbings.

\para{Donaldson's theorem.}
We applied the criterion of Jabuka and Kelly \cite[Thms.~2.10 and 2.12]{JK} for the Goeritz forms of alternating
diagrams.  Among the $\nJKTargets$
alternating knots with at most $13$ crossings whose range starts at $1$, it excludes exactly the $\nJKObstr$ knots
excluded by the linking form ($\nJKLinking$) or by Proposition~\ref{prop:dtest} ($\nDAltObstr$); among the
$\nLSJKTargets$ alternating knots with $14$ crossings left open by \cite{LS}, it excludes exactly the $\nLSJKObstr$
knots excluded by the linking form ($\nLSJKLinking$) or by Proposition~\ref{prop:dtest} ($\nLSDObstr$).  So on these
knots the correction terms give no more than the linking form and Donaldson's theorem, as \S\ref{sec:donaldson}
predicts.  The method of Proposition~\ref{prop:newd}(c) was also applied to the $\nDefControls$ non-alternating
knots that are not Montesinos, have $\sigma+4\Arf\equiv\pm2\pmod 8$ and bound the smooth M\"obius bands of
Proposition~\ref{prop:upper}(a): for $\nDefControlsSurface$ of them some diagram has a
definite surface of the required sign, and all $\nDefControlLattices$ of these lattices embed.  The sign matters: for
the Goeritz forms of alternating diagrams, the criterion excludes none of the $\nJKControls$ alternating knots with \KI{} value $1$ above, but with the opposite sign it would exclude $\nJKReversed$ of them.

\subsection{Knots with \texorpdfstring{$h_\lambda=3$}{h = 3} and the rational Witt span}\label{sec:ws}
Jabuka's rational Witt span $\ws(K)$ is the least rank of a form over $\Q$ whose residues at all primes agree with those
of the rational Witt class of $V+V^{T}$ \cite[Def.~1.3]{Jabuka}; it bounds the concordance crosscap number from below
\cite[Thm.~1.4]{Jabuka}.  By \cite[Thm.~4.6]{LeeTheory}, $\ws(K)\le h_\lambda(K)\le\ws(K)+1$.  Moreover
$\ws(K)=3$ if and only if (i) or (ii) of Corollary~\ref{cor:criteria}(c) holds, or some $p\equiv3\pmod 4$ with $a_p=2$
has $\varepsilon_p=-1$ and $(\sigma(K),D)\equiv(4,1)\pmod 8$; and $h_\lambda(K)=\ws(K)+1$ if and only if (i) and (ii)
fail, some $p\equiv3\pmod4$ with $a_p=2$ has $\varepsilon_p=-1$, and $(\sigma(K),D)\equiv(0,5)\pmod 8$.

In the table of the $\nHthree$ knots with at most $13$ crossings and $h_\lambda=3$ (Table~\ref{tab:hthree}), whose
column $\ws$ gives these values, the last case applies exactly to the $\nWsTwo$ knots that satisfy only (iii) and have
$\sigma(K)\equiv0\pmod 8$, among them $\kn{8}{18}$; these knots have $\ws=2$.  The other $\nWsThree$ knots of the table
have $\ws=3$: $\nWsThreeOne$ of them satisfy (i), and the $\nWsThreeIII$ knots $\kn{12a}{569}$, $\kn{12n}{626}$,
$\kn{13a}{1974}$, $\kn{13a}{3121}$, $\kn{13n}{3005}$ and $\kn{13n}{4927}$ satisfy (iii) with $\sigma(K)\equiv4\pmod 8$.
So by \cite[Thm.~4.6(a)]{LeeTheory}, for these six knots no binary form $Q$ of odd determinant, odd or even,
has $[\lambda_Q]=[\lambda_K]$.  No knot of the table satisfies (ii).  Hence $\ws(K)=h_\lambda(K)$ for all but $\nWsTwo$
of the $\nKnotsThirteen$ prime knots with at most $13$ crossings.  \cite[Example~1.7]{Jabuka} lists the knots with at
most $12$ crossings to which \cite[Thm.~1.6(ii)]{Jabuka} applies; that criterion requires the residues of the rational
Witt class to vanish at all primes but one.  They are exactly the knots of the table with at most $12$ crossings that
satisfy (i).  The knots $\kn{12a}{569}$ and $\kn{12n}{626}$ also have $\ws=3$, although the criterion does not apply to
them.  Of the $\nLSHthree$ knots with $14$ crossings and $h_\lambda=3$ (Theorem~\ref{mthm:fourteen}), $\nLSWsTwo$
have $\ws=2$, and the others have $\ws=3$.  A separate computation of $\ws$ from the residues alone, which decides the
existence of rational forms of rank $2$ with Hilbert symbols and the Hasse--Minkowski theorem, gives the same values.

\subsection{Fourteen crossings}\label{sec:fourteen}
For the knots with $14$ crossings we give only lower bounds, from the linking form and, for the alternating knots, from
Proposition~\ref{prop:dtest}.

\begin{proof}[Proof of Theorem~\ref{mthm:fourteen}]
Of the $\nLSOpen$ knots left open by \cite{LS}, $\nLSOpenLinking$ have $h_\lambda\ge2$ and hence $\gfour\ge2$ by
Theorem~\ref{thm:lower}, and among the others Proposition~\ref{prop:dtest}, with the correction terms of Goeritz
manifolds (\S\ref{sec:dtest}), shows $\gfour\ge2$ for $\nLSDObstr$ of the $\nLSDTargets$ alternating knots; so
$\gfour\ge2$ for $\nLSNew$ of them.  The $\nLSHthree$ knots with $h_\lambda=3$, $\nLSHthreeAlt$ of them alternating,
have $\gtop\ge3$ by Theorem~\ref{thm:lower}.
\end{proof}

The linking forms were computed from Goeritz matrices of the SnapPy diagrams, which are alternating for the
alternating knots.  They agree with the table of \cite{LS}: none of the $\nLSBand$ knots for which Lee and Sabloff \cite{LS} found a M\"obius band has $h_\lambda\ge2$, and every
knot that they obstruct by Yasuhara's condition ($\nLSArf$ knots) or by the linking form ($\nLSLinking$ knots) has
$h_\lambda\ge2$.  The correction-term test excludes none of the $\nLSDControls$ alternating knots with a M\"obius band.

\subsection{Corrections to \texorpdfstring{\KI}{KnotInfo}}\label{sec:corrections}

For five knots with $11$ or fewer crossings our results differ from \KI{} 2026.11.2, from the sources it draws on, or
from the arguments it cites.

\para{$\kn{9}{49}$.}  \KI{} lists $\gfour=2$ and $\gtop\in[1,2]$.  Here $H_K\cong\Z/5\oplus\Z/5$ and condition~(i) of
Corollary~\ref{cor:criteria}(c) holds, which in this case is \cite[Thm.~4]{GL}; so $\gtop\ge3$.  Jabuka and Kelly show
$\gfour(\kn{9}{49})=3$ in the revised version \cite[Thm.~1.2 and \S3.2]{JKv3}; the published version
\cite{JK} lists $2$, from a band move to $\kn{8}{21}$ and the value $\gfour(\kn{8}{21})=1$, which their erratum
\cite{JKerr} corrects to $2$.  Since the crosscap number of $\kn{9}{49}$ is $3$ \cite{KnotInfo}, $\gfour=\gtop=3$.

\para{$\kn{11n}{133}$.}  \KI{} lists $\gfour=2$ and $\gtop\in[1,2]$, from \cite[Prop.~3.3 and Lemma~3.7]{Fairchild}:
the lower bound is Yasuhara's condition $\sigma(K)+4\Arf(K)\equiv4\pmod 8$ \cite[Prop.~2.3]{Fairchild}, and the upper
bound is a nonorientable band move to $\kn{10}{165}$ \cite[Fig.~17(d)]{Fairchild}, for which \KI{} lists the value $1$.
But $H_K\cong\Z/5\oplus\Z/5$ and Corollary~\ref{cor:criteria}(c)(i) applies, so $\gtop\ge3$, and $\kn{11n}{133}$ has no
band move to a knot that bounds a M\"obius band, such as $\kn{10}{165}$, which a band move turns into the unknot
\cite[Thm.~1.1 and Fig.~35(i)]{Ghanbarian}; so the band move of \cite[Fig.~17(d)]{Fairchild} gives a different knot.  In
the \KI{} diagram of $\kn{11n}{133}$ a band with $t=0$ gives the figure-eight knot $\kn{4}{1}$, and two further bands
turn it into the unknot, as for $\kn{13a}{2890}$ in \S\ref{sec:upper}.  With a disc, the three bands form a
nonorientable surface with $b_1=3$, so $\gfour(\kn{11n}{133})=\gtop(\kn{11n}{133})=3$.  In $\nElevenControlDiagrams$
diagrams of $\kn{11n}{133}$ we found no band to a knot with \KI{} value $1$.

\para{$\kn{11n}{176}$.}  \KI{} lists $\gfour=2$, with the lower bound from \cite[Lemma~3.9]{Fairchild}.  That lemma
concerns $H_K\cong\Z/p^2q$ with every prime of $q$ occurring to an odd power, and its conclusion is that some generator $a$ has
$\lambda(a,a)=\pm1/p^2q$ or $\lambda(a,a)=\pm1/q$.  The second alternative is impossible, since $\lambda(a,a)$ has
denominator $p^2q$; in the case where the $p$-primary part lies in the metabolic summand, the proof shows only that the
form on the summand $G_1$ of the (Murakami--Yasuhara) splitting is $\pm1/q$.  For $\kn{11n}{176}$, with
$H_K\cong\Z/63$ and $\lambda(x,x)=2/63$ for a generator $x$, the first alternative fails, as neither $2$ nor $-2$ is a
square modulo $63$, but this weaker conclusion holds.  In fact no argument that uses only the linking form, the
signature modulo $8$ and the Arf invariant can give $\gtop\ge2$ here: the connected sum of $\kn{7}{1}$ and the mirror
image of $\kn{6}{1}$ bounds a M\"obius band and has a linking form isometric to $\lambda_{\kn{11n}{176}}$, the same
signature modulo $8$ (its signature is $-6$, and that of $\kn{11n}{176}$ is $2$), and the same Arf invariant.  Indeed,
$\lambda(9x,9x)=4/7$ and $\lambda(7x,7x)=5/9$, which are values on generators of $\lambda_{\kn{7}{1}}$ and of
$-\lambda_{\kn{6}{1}}$, the form of the mirror image of $\kn{6}{1}$.

The correction terms of $\Sigma_K$ do not decide the question either.  The knot $\kn{11n}{176}$ has unknotting number
one, and the method of Proposition~\ref{prop:newd}(b) gives $\Sigma_K\cong S^3_{63/2}(C)$, where $C$ is a fibred knot of genus $9$
with $\tau(C)=9$ whose exterior is the manifold \textup{K6\_29} in the tables of SnapPy.  The reduced Khovanov homology of
$\kn{11n}{176}$ over $\Z/2$, as listed by \KI{}, has rank $63=D$, so $\Sigma_K$ is an $L$-space
\cite[Cor.~1.2]{OSz-bdc}.  Hence $C$ has a positive $L$-space surgery, the correction terms of $S^3_{63/2}(C)$ are given
by the torsion coefficients of the Alexander polynomial of $C$ \cite[Thm.~1.2]{OSz-rat}, and
\[
(V_0,\dots,V_9)(C)=(4,3,3,2,2,2,1,1,1,0),\qquad V_j(C)=0\quad(j\ge9).
\]
The meridian $\mu$ of $C$ generates $H_K$, and $\lambda_K(\mu,\mu)=-\lambda_{\Sigma_K}(\mu,\mu)=2/63$
(\S\ref{sec:dtest}).  Take $Z=-\Sigma_K$, $\ell=7$, $M=\langle21\mu\rangle$ and $g=6\mu$ in
Proposition~\ref{prop:dtest}: $M$ is isotropic, $g$ generates $M^\perp/M=\langle3\mu\rangle/M\cong\Z/7$, and
$\lambda_Z(g,g)=\lambda_K(6\mu,6\mu)=8/7\equiv1/7$.  The formula of Ni and Wu (\S\ref{sec:dtest}) gives, for every
$m\in M$,
\[
\begin{array}{c|ccccccc}
j&0&1&2&3&4&5&6\\\hline
4d(Z,\mathfrak t_s+jg+m)&2&-18/7&-2/7&6/7&6/7&-2/7&-18/7\\
1-(2j-7)^2/7&-6&-18/7&-2/7&6/7&6/7&-2/7&-18/7
\end{array}
\]
Thus \eqref{eq:dtest} holds, with equality for $j\ge1$.  It fails for every other choice, since $\ell|M|^2=63$ leaves
$(\ell,|M|)=(63,1)$ or $(7,3)$.  For $\ell=63$ a generator $g=a\mu$ of $H_K$ would need $\pm2a^2\equiv1\pmod{63}$, with
the sign $+$ for $Z=-\Sigma_K$, which fails modulo $3$, and $-$ for $Z=\Sigma_K$, which fails modulo $7$.  For $\ell=7$ the
subgroup $M=\langle21\mu\rangle$ is the only one of order $3$, and $g=3a\mu$ needs $\pm2a^2\equiv1\pmod 7$, which fails
for $Z=\Sigma_K$ and gives $g\equiv\pm6\mu$ modulo $M$ for $Z=-\Sigma_K$.  Replacing $g$ by $-g$ or by $g+m$ with $m\in M$
gives the same inequalities, as $d(Z,\mathfrak t_s+x)=d(Z,\mathfrak t_s-x)$ \cite[Prop.~4.2]{OSz-abs} and $m$ runs over
$M$.
So $\gfour(\kn{11n}{176})\in[1,2]$ remains open; the double branched cover $W$ of a smooth M\"obius band would have
$Q_W=(7)$, with $7>0$ by the proof of Proposition~\ref{prop:donaldson} since $\sigma(K)+4\Arf(K)\equiv2\pmod 8$, and
$3$ dividing $|H_1(W)|$.

\para{$\kn{11n}{38}$.}  \KI{} lists $\gfour=2$, from \cite[Lemma~4.1]{Fairchild}.  Its proof takes the double branched
cover of a M\"obius band bounded by $\kn{11n}{38}$ to be negative definite, citing \cite{GL}, and then infers from the
linking form, read off from a Goeritz matrix, that the cover is positive definite.  The trefoil $\kn{3}{1}$ of \KI{} has
the same linking form $\langle1/3\rangle$, signature $-2$ and Arf invariant as $\kn{11n}{38}$ and bounds a M\"obius
band, so these invariants alone cannot decide the question.  Proposition~\ref{prop:dtest} does: $\kn{11n}{38}$ is the
Montesinos knot $K(1/2;-2/3;2/9)$, and its correction terms exclude a smooth M\"obius band, so
$\gfour(\kn{11n}{38})\ge2$.  The upper bound $2$ rests on the band move to $\kn{3}{1}$ of \cite[Lemma~4.1]{Fairchild},
which we have not rebuilt; $\gtop(\kn{11n}{38})\in[1,2]$.

\para{$\kn{11n}{57}$.}  \KI{} 2026.11.2 lists the range $[1,3]$.  But \cite[Prop.~3.2]{Fairchild} gives
$\gfour(\kn{11n}{57})=1$, by a band move with $t=0$ to the unknot \cite[Fig.~8(g)]{Fairchild}, and so do the table of \cite{LS} and version 2026.8.1 of \KI{}; Lee and Sabloff found that this band move gives the figure-eight
knot, and they leave the value open \cite[\S4.3.2]{LS}. Here $H_K\cong\Z/7$ and $h_\lambda(K)=1$, and $\kn{11n}{57}$ is
the Montesinos knot $K(1/2;-2/3;3/11)$. Its correction terms exclude a smooth M\"obius band
(Proposition~\ref{prop:dtest}), both with the plumbing given by this notation and with the star-shaped plumbing built from the Seifert
invariants that Regina recognises (\S\ref{sec:gap}).  With the punctured Klein bottle
of Proposition~\ref{prop:upper}(b), $\gfour(\kn{11n}{57})=2$, while $\gtop(\kn{11n}{57})\in[1,2]$.

\section{Smooth versus locally flat M\"obius bands}\label{sec:gap}

Theorem~\ref{thm:lower} concerns locally flat surfaces, while Proposition~\ref{prop:dtest} uses smoothness.  A knot with
$h_\lambda=1$ that fails the test of Proposition~\ref{prop:dtest} may therefore bound a locally flat M\"obius band
without bounding a smooth one.  By Proposition~\ref{prop:lower}(c), Proposition~\ref{prop:dtest} excludes a
smooth M\"obius band for $\nGapCandidates$ knots with \KI{} range $[1,2]$ ($\nGapAlt$ alternating and $\nGapMont$
non-alternating Montesinos knots), all with $h_\lambda=1$.  We searched for locally flat M\"obius bands for these $\nGapCandidates$ knots and for all $\nNewD$ knots of Proposition~\ref{prop:newd}(a), (b), which also have $h_\lambda=1$; we did not test the $\nDObstrOneThree$ knots with range $[1,3]$ excluded by Proposition~\ref{prop:lower}(c) or the $\nDefCert$ knots of Proposition~\ref{prop:newd}(c).

Twelve of the $\nGapAll$ knots of Theorem~\ref{mthm:gap} are the knots
\begin{gather*}
\kn{12a}{820},\ \kn{12n}{25},\ \kn{12n}{198},\ \kn{12n}{550},\ \kn{13a}{26},\ \kn{13a}{938},\\
\kn{13a}{945},\ \kn{13a}{1601},\ \kn{13a}{3543},\ \kn{13a}{4039},\ \kn{13a}{4127},\ \kn{13a}{4578}
\end{gather*}
of Proposition~\ref{prop:lower}(c), and the other $\nGapNew$ are knots of Proposition~\ref{prop:newd}(a), (b) with $13$
crossings.

\begin{proof}[Proof of Theorem~\ref{mthm:gap}]
For each $K$, $h_\lambda(K)=1$, so $\gtop(K)\ge1$ by Theorem~\ref{thm:lower}.  A band move turns $K$ into a knot $J$
with Alexander polynomial $\Delta_J=1$ (Table~\ref{tab:gap}).  As $J$ is a knot, the band is nonorientable
(\S\ref{sec:bandconv}).  By Freedman's theorem \cite[Thm.~11.7B and Cor.~7.1B]{FQ}, $J$ bounds a locally flat disc $E$ in
$B^4$.  The band, pushed into a
collar, together with $E$ is a locally flat nonorientable surface of Euler characteristic $0$ with boundary $K$, that
is, a M\"obius band.  So $\gtop(K)=1$.  Finally $K$ bounds no smooth M\"obius band by Proposition~\ref{prop:dtest}:
\eqref{eq:dtest} fails for both choices of $Z$ and for every $\ell$, $M$ and $g$, by
Propositions~\ref{prop:lower}(c) and~\ref{prop:newd}(a), (b) (for the knots of Proposition~\ref{prop:newd}(b), for every
$V$ between the bounds of \S\ref{sec:dtest}, or for the exact $V_j(C)$).  So $\gfour(K)\ge2$, and with the upper bound
$2$ of \KI{} or of Proposition~\ref{prop:upper}(b), $\gfour(K)=2$.
\end{proof}

\def\gapcaption{The knots of Theorem~\ref{mthm:gap}.  The column ``$d$ from'' names the source of the correction terms:
the Goeritz manifold or the star-shaped plumbing of Proposition~\ref{prop:lower}(c), or the plumbing tree or the surgery
description of Proposition~\ref{prop:newd}(a), (b).  A knot $J$ is named as in \KI{} or SnapPy when SnapPy identifies it,
and otherwise by its hyperbolic volume as computed numerically by SnapPy, which serves only as a label.  The proof uses
one band for each $K$, and its knot $J$ is given by a PD code in the repository.}
\begin{table}[p]
\centering
\footnotesize
\caption{\gapcaption}\label{tab:gap}
\begin{tabular}{lrll}
\toprule
$K$ & $\det K$ & $d$ from & knots $J$ with $\Delta_J=1$ one band move away \\
\midrule
$\kn{12a}{820}$ & $179$ & Goeritz & volume $19.4549$ \\
$\kn{12n}{25}$ & $11$ & star plumbing & $\kn{13n}{1496}$; volume $17.4178$; \textup{K14n7708} \\
$\kn{12n}{198}$ & $27$ & star plumbing & volume $11.8463$ \\
$\kn{12n}{550}$ & $33$ & star plumbing & $\kn{13n}{1496}$ \\
$\kn{13a}{26}$ & $209$ & Goeritz & $\kn{13n}{1496}$ \\
$\kn{13a}{938}$ & $275$ & Goeritz & volume $17.6576$; volume $23.0126$ \\
$\kn{13a}{945}$ & $275$ & Goeritz & $\kn{13n}{1496}$; \textup{K14n7708} \\
$\kn{13a}{1601}$ & $187$ & Goeritz & $\kn{13n}{1496}$ \\
$\kn{13a}{3543}$ & $331$ & Goeritz & volume $23.6119$ \\
$\kn{13a}{4039}$ & $267$ & Goeritz & volume $16.8017$ \\
$\kn{13a}{4127}$ & $283$ & Goeritz & $\kn{13n}{1496}$ \\
$\kn{13a}{4578}$ & $127$ & Goeritz & volume $15.9543$ \\
$\kn{13n}{472}$ & $11$ & surgery & \textup{K14n7708}; volume $22.3337$ \\
$\kn{13n}{594}$ & $87$ & plumbing tree & $\kn{13n}{1496}$ \\
$\kn{13n}{628}$ & $3$ & plumbing tree & $\kn{13n}{1496}$; volume $17.1857$ \\
$\kn{13n}{1003}$ & $39$ & plumbing tree & $\kn{13n}{1496}$ \\
$\kn{13n}{1134}$ & $11$ & plumbing tree & volume $16.8017$; volume $20.4976$ \\
$\kn{13n}{1164}$ & $7$ & plumbing tree & volume $15.9543$ \\
$\kn{13n}{1281}$ & $3$ & plumbing tree & volume $18.8415$ \\
$\kn{13n}{1307}$ & $3$ & plumbing tree & volume $18.8415$ \\
$\kn{13n}{1411}$ & $31$ & plumbing tree & \textup{K14n7708} \\
$\kn{13n}{1437}$ & $31$ & plumbing tree & volume $17.5308$ \\
$\kn{13n}{1576}$ & $3$ & plumbing tree & volume $17.6576$; $\kn{13n}{1496}$ \\
$\kn{13n}{1618}$ & $11$ & plumbing tree & volume $22.2606$; $\kn{13n}{1496}$ \\
$\kn{13n}{1861}$ & $3$ & surgery & volume $15.4897$; volume $13.2191$; volume $17.9280$ \\
$\kn{13n}{2010}$ & $57$ & surgery & volume $23.9682$ \\
$\kn{13n}{2032}$ & $3$ & surgery & volume $16.6453$ \\
$\kn{13n}{2182}$ & $3$ & plumbing tree & volume $13.2263$ \\
$\kn{13n}{2206}$ & $27$ & surgery & volume $21.4638$ \\
$\kn{13n}{2271}$ & $7$ & plumbing tree & volume $20.0050$ \\
$\kn{13n}{2352}$ & $83$ & plumbing tree & $\kn{13n}{1496}$ \\
$\kn{13n}{2366}$ & $47$ & plumbing tree & \textup{K14n7708} \\
$\kn{13n}{2386}$ & $59$ & plumbing tree & $\kn{13n}{1496}$ \\
$\kn{13n}{2932}$ & $51$ & plumbing tree & volume $20.0163$ \\
$\kn{13n}{3183}$ & $43$ & surgery & \textup{K14n7708} \\
$\kn{13n}{3234}$ & $11$ & surgery & volume $21.2019$; volume $19.2660$ \\
$\kn{13n}{3249}$ & $31$ & plumbing tree & volume $11.8463$ \\
$\kn{13n}{3547}$ & $11$ & plumbing tree & volume $14.3908$ \\
$\kn{13n}{3844}$ & $11$ & surgery & volume $24.3191$; $\kn{13n}{1496}$ \\
$\kn{13n}{3855}$ & $11$ & surgery & volume $24.3191$ \\
$\kn{13n}{4947}$ & $151$ & surgery & volume $19.9314$ \\
\bottomrule
\end{tabular}
\end{table}

\begin{figure}[ht]
\centering
\begin{minipage}[b]{0.48\textwidth}\centering
\includegraphics{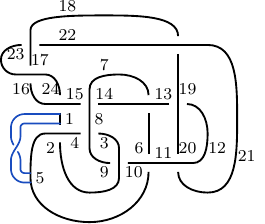}\\[2pt]
{\small $\kn{12n}{25}$}
\end{minipage}\hfill
\begin{minipage}[b]{0.48\textwidth}\centering
\includegraphics{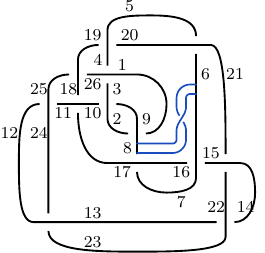}\\[2pt]
{\small $\kn{13a}{1601}$}
\end{minipage}

\bigskip
\begin{minipage}[b]{0.48\textwidth}\centering
\includegraphics{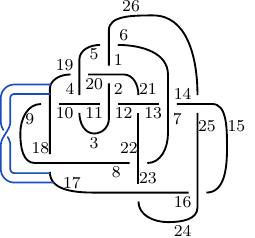}\\[2pt]
{\small $\kn{13a}{4127}$}
\end{minipage}\hfill
\begin{minipage}[b]{0.48\textwidth}\centering
\includegraphics{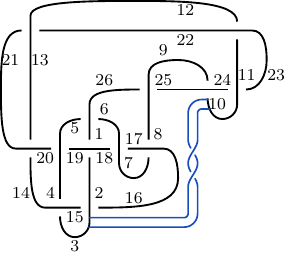}\\[2pt]
{\small $\kn{13a}{4578}$}
\end{minipage}
\caption{The bands of Theorem~\ref{mthm:gap} in the \KI{} diagrams of four of the twelve knots of
Proposition~\ref{prop:lower}(c), drawn from their PD codes
with the edges numbered; the face, the two edges and $t$ (as in Figure~\ref{fig:band}) of each band are listed in
\S\ref{sec:gap}.  The band move
along each band gives a knot with Alexander polynomial one.}
\label{fig:gap}
\end{figure}

For four of the twelve knots of Proposition~\ref{prop:lower}(c) a band in the \KI{} diagram itself works.  In the notation
of Figure~\ref{fig:band} (edges of the \KI{} PD code, and $t$ the signed number of crossings between the long sides of
the band) these bands are as follows.
\begin{itemize}
\item $\kn{12n}{25}$: the face with edges $1, 16, 23, 18, 21, 11, 5$, joining edges $1$ and $5$, $t=-1$, gives $\kn{13n}{1496}$.
\item $\kn{13a}{1601}$: the face with edges $1, 8, 16, 6, 20, 4$, joining edges $8$ and $6$, $t=-1$, gives $\kn{13n}{1496}$.
\item $\kn{13a}{4127}$: the face with edges $26, 15, 24, 17, 9, 19, 5$, joining edges $17$ and $19$, $t=-1$, gives $\kn{13n}{1496}$.
\item $\kn{13a}{4578}$: the face with edges $3, 14, 21, 12, 23, 10, 25, 8, 16$, joining edges $3$ and $10$, $t=-2$, gives the knot of numerically computed volume $15.9543$.
\end{itemize}
These bands are drawn in Figure~\ref{fig:gap}.  For $\nGapNewKI$ of the $\nGapNew$ knots of Proposition~\ref{prop:newd}
the band lies in the \KI{} diagram as well; all bands are recorded in the repository.

\para{The bands.}
We looked for such bands in $\nGapDiagrams$ diagrams of these $\nGapCandidates$ knots: the \KI{} diagram, diagrams found by
SnapPy through random Reidemeister moves and simplification, and further diagrams obtained from these by random
Reidemeister moves, all with at most six crossings more than the crossing number and each verified to have the
complement of the \KI{} knot.  In every face we tried the bands between any two of its edges
with $|t|\le6$.  For a fixed arc of the band, $\det J_t=|\alpha+t\det K|$ for an integer $\alpha$, so, as $\det K\ge3$, at most
one $t$ gives $\det J_t=1$; when it lay outside this range we tried it as well.  Of the $\nGapBands$ bands, $\nGapKnots$ gave
knots, $\nGapDetOne$ of them with determinant $1$ and $\nGapHits$ with $\Delta_J=1$.  The Alexander polynomials were
computed by Fox calculus and verified with a Seifert matrix and with knot Floer homology.  Each of these $J$ has one
component and gives back $K$ by the dual band.  The twelve knots of Theorem~\ref{mthm:gap} from
Proposition~\ref{prop:lower}(c) are those with such a band.  Applied to the $\nGapNewCandidates$ knots of
Proposition~\ref{prop:newd}, the same procedure found bands for $\nGapNew$ of them ($\nGapNewPlumb$ from (a) and $\nGapNewSurg$ from (b)
only), and all resulting knots $J$ with $\Delta_J=1$ were checked in the same way.

For $\nGapNoForm$ of the $\nGapCandidates$ knots no such band exists.  By the Montesinos trick \cite{Montesinos}, a band move between
$K$ and $J$ presents $\Sigma_K$ as integral surgery on a knot in $\Sigma_J$, and if $\det J=1$ this forces
$\lambda_K\cong\langle\pm1/D\rangle$.  For these $\nGapNoForm$ knots $\lambda_K$ is not of this form, and indeed no knot
$J$ with $\det J=1$ occurred.  For $\nGapRestDetOne$ of the remaining $\nGapRest$ knots we found knots $J$
with $\det J=1$ but none with $\Delta_J=1$, and for the other $\nGapRestNoDetOne$ we found no $J$ with $\det J=1$.  These diagrams do not exhaust the possibilities, and larger diagrams may give bands for more of them.

\para{Independent verification.}
For each of the twelve knots, one band was reconstructed in a separate check: the diagram was confirmed to represent $K$, the
band move was repeated both directly and via the band operation of SnapPy's link module, and $\Delta_J=1$ was confirmed
from a Seifert matrix and from knot Floer homology.  Every $J$ that occurred for these twelve knots has
Ozsv\'ath--Szab\'o invariant $\tau(J)=\pm1$, so none of them is smoothly slice \cite{OSz-tau}, in accordance with
Theorem~\ref{mthm:gap}.  For the $\nGapNew$ knots of Proposition~\ref{prop:newd} we recomputed, for the band used in the
proof, the knot Floer homology of $J$: $\Delta_J=1$, and $\tau(J)=\pm1$ except for $\kn{13n}{1281}$, $\kn{13n}{1307}$,
$\kn{13n}{2010}$ and $\kn{13n}{2932}$, whose knots $J$ have $\tau(J)=0$.  These $\nGapTauZeroJ$ knots $J$
are pairwise distinct: the total ranks of their knot Floer homology are $97$, $97$, $257$ and $145$, and the two of
rank $97$, from $\kn{13n}{1281}$ and $\kn{13n}{1307}$, have Seifert genera $3$ and $4$.  They are topologically slice, and by Theorem~\ref{mthm:gap} they are
not smoothly slice, since a smooth slice disc of $J$ and the band would form a smooth M\"obius band with boundary $K$.
The smooth obstructions were also checked by a second route.  For the nine alternating knots of
Proposition~\ref{prop:lower}(c), the lattice criterion of Jabuka and Kelly gives the same conclusion
(\S\ref{sec:controls}).  For the three Montesinos knots, $\Sigma_K$ was identified as a Seifert fibred space with Regina
\cite{Regina}, and its plumbing and correction terms were recomputed from scratch; this calculation excludes all
three, excludes none of the $\nMontCheckControls$ Montesinos knots with \KI{} value $1$ that it can treat, and
reproduces $d(\Sigma_K,\mathfrak t_s)$ from the first computation for all of them.  For the knots of
Proposition~\ref{prop:newd} the verification is described in \S\ref{sec:lower}.

\subsection*{Data availability and AI declaration.}
The codes and tables for the computations in this paper, including the explicit band moves, plumbing trees,
surgery descriptions and definite surfaces, are available at
\begin{center}
\coderepo
\end{center}
The computations use SnapPy and its
link module spherogram \cite{SnapPy}, Regina \cite{Regina}, PARI \cite{PARI}, SageMath \cite{Sage} and Szab\'o's program
for knot Floer homology \cite{Szabo}.

The AI model Claude Opus 5.5 (Anthropic) was partially used in the
preparation of this paper to write the code for the computations and to check the grammar of the draft.

\subsection*{Acknowledgements}
This work is supported by the Institute for Basic Science (IBS) under Project No.\ IBS-R003-D1.

\FloatBarrier
\appendix

\Needspace{10\baselineskip}
\section{Tables}\label{app:tables}

\subsection{Knots with \texorpdfstring{$\gfour=2$}{gamma4 = 2} and range \texorpdfstring{$[1,2]$}{[1,2]}}\label{app:newtwo}
The following $\nGammaTwoOneTwo$ knots, whose \KI{} range is $[1,2]$, have $\gfour=2$ (\S\ref{sec:values}).  The entries are ordered down the columns.

{\footnotesize
\setlength{\columnsep}{4pt}
\setlength{\multicolsep}{3pt}
\par\vspace{3pt}\hrule height\heavyrulewidth
\begin{multicols}{7}
\raggedright\setlength{\parindent}{0pt}\setlength{\parskip}{0pt}
$\kn{11a}{6}$\par
$\kn{11a}{59}$$^{\mathrm{d}}$\par
$\kn{11a}{173}$\par
$\kn{11a}{188}$$^{\mathrm{d}}$\par
$\kn{11a}{264}$\par
$\kn{11a}{277}$\par
$\kn{11a}{280}$\par
$\kn{11a}{352}$\par
$\kn{12a}{389}$\par
$\kn{12a}{426}$\par
$\kn{12a}{503}$\par
$\kn{12a}{551}$$^{\mathrm{d}}$\par
$\kn{12a}{595}$$^{\mathrm{d}}$\par
$\kn{12a}{712}$\par
$\kn{12a}{820}$$^{\mathrm{d,t}}$\par
$\kn{12a}{969}$$^{\mathrm{d}}$\par
$\kn{12a}{981}$\par
$\kn{12a}{1015}$$^{\mathrm{d}}$\par
$\kn{12a}{1242}$$^{\mathrm{d}}$\par
$\kn{12a}{1256}$$^{\mathrm{d}}$\par
$\kn{12a}{1276}$$^{\mathrm{d}}$\par
$\kn{12a}{1279}$$^{\mathrm{d}}$\par
$\kn{12n}{25}$$^{\mathrm{d,t}}$\par
$\kn{12n}{41}$$^{\mathrm{d}}$\par
$\kn{12n}{47}$$^{\mathrm{d}}$\par
$\kn{12n}{64}$$^{\mathrm{d}}$\par
$\kn{12n}{79}$$^{\mathrm{d}}$\par
$\kn{12n}{198}$$^{\mathrm{d,t}}$\par
$\kn{12n}{249}$$^{\mathrm{d}}$\par
$\kn{12n}{261}$$^{\mathrm{d}}$\par
$\kn{12n}{550}$$^{\mathrm{d,t}}$\par
$\kn{12n}{637}$\par
$\kn{12n}{799}$\par
$\kn{12n}{869}$\par
$\kn{13a}{26}$$^{\mathrm{d,t}}$\par
$\kn{13a}{32}$$^{\mathrm{d}}$\par
$\kn{13a}{33}$$^{\mathrm{d}}$\par
$\kn{13a}{58}$\par
$\kn{13a}{140}$$^{\mathrm{d}}$\par
$\kn{13a}{189}$$^{\mathrm{d}}$\par
$\kn{13a}{223}$$^{\mathrm{d}}$\par
$\kn{13a}{316}$\par
$\kn{13a}{319}$$^{\mathrm{d}}$\par
$\kn{13a}{349}$\par
$\kn{13a}{361}$\par
$\kn{13a}{395}$\par
$\kn{13a}{428}$$^{\mathrm{d}}$\par
$\kn{13a}{445}$\par
$\kn{13a}{553}$$^{\mathrm{d}}$\par
$\kn{13a}{635}$\par
$\kn{13a}{654}$\par
$\kn{13a}{662}$\par
$\kn{13a}{663}$$^{\mathrm{d}}$\par
$\kn{13a}{675}$$^{\mathrm{d}}$\par
$\kn{13a}{676}$$^{\mathrm{d}}$\par
$\kn{13a}{702}$$^{\mathrm{d}}$\par
$\kn{13a}{738}$\par
$\kn{13a}{740}$\par
$\kn{13a}{762}$\par
$\kn{13a}{765}$\par
$\kn{13a}{813}$$^{\mathrm{d}}$\par
$\kn{13a}{831}$$^{\mathrm{d}}$\par
$\kn{13a}{832}$\par
$\kn{13a}{892}$$^{\mathrm{d}}$\par
$\kn{13a}{938}$$^{\mathrm{d,t}}$\par
$\kn{13a}{945}$$^{\mathrm{d,t}}$\par
$\kn{13a}{968}$\par
$\kn{13a}{971}$\par
$\kn{13a}{993}$\par
$\kn{13a}{1022}$$^{\mathrm{d}}$\par
$\kn{13a}{1089}$$^{\mathrm{d}}$\par
$\kn{13a}{1100}$\par
$\kn{13a}{1221}$$^{\mathrm{d}}$\par
$\kn{13a}{1226}$$^{\mathrm{d}}$\par
$\kn{13a}{1266}$\par
$\kn{13a}{1395}$\par
$\kn{13a}{1429}$$^{\mathrm{d}}$\par
$\kn{13a}{1474}$$^{\mathrm{d}}$\par
$\kn{13a}{1488}$$^{\mathrm{d}}$\par
$\kn{13a}{1573}$$^{\mathrm{d}}$\par
$\kn{13a}{1601}$$^{\mathrm{d,t}}$\par
$\kn{13a}{1640}$$^{\mathrm{d}}$\par
$\kn{13a}{1682}$\par
$\kn{13a}{1720}$\par
$\kn{13a}{1751}$$^{\mathrm{d}}$\par
$\kn{13a}{1786}$$^{\mathrm{d}}$\par
$\kn{13a}{1791}$$^{\mathrm{d}}$\par
$\kn{13a}{1829}$\par
$\kn{13a}{1855}$$^{\mathrm{d}}$\par
$\kn{13a}{1879}$$^{\mathrm{d}}$\par
$\kn{13a}{1897}$\par
$\kn{13a}{1899}$\par
$\kn{13a}{1971}$\par
$\kn{13a}{1976}$$^{\mathrm{d}}$\par
$\kn{13a}{2008}$\par
$\kn{13a}{2109}$$^{\mathrm{d}}$\par
$\kn{13a}{2114}$$^{\mathrm{d}}$\par
$\kn{13a}{2202}$\par
$\kn{13a}{2207}$\par
$\kn{13a}{2307}$\par
$\kn{13a}{2315}$\par
$\kn{13a}{2367}$$^{\mathrm{d}}$\par
$\kn{13a}{2381}$\par
$\kn{13a}{2432}$$^{\mathrm{d}}$\par
$\kn{13a}{2461}$$^{\mathrm{d}}$\par
$\kn{13a}{2484}$\par
$\kn{13a}{2497}$$^{\mathrm{d}}$\par
$\kn{13a}{2501}$\par
$\kn{13a}{2544}$\par
$\kn{13a}{2672}$\par
$\kn{13a}{2690}$$^{\mathrm{d}}$\par
$\kn{13a}{2717}$\par
$\kn{13a}{2722}$$^{\mathrm{d}}$\par
$\kn{13a}{2736}$$^{\mathrm{d}}$\par
$\kn{13a}{2742}$\par
$\kn{13a}{2933}$\par
$\kn{13a}{2947}$\par
$\kn{13a}{3016}$$^{\mathrm{d}}$\par
$\kn{13a}{3022}$$^{\mathrm{d}}$\par
$\kn{13a}{3036}$$^{\mathrm{d}}$\par
$\kn{13a}{3131}$\par
$\kn{13a}{3143}$$^{\mathrm{d}}$\par
$\kn{13a}{3474}$$^{\mathrm{d}}$\par
$\kn{13a}{3543}$$^{\mathrm{d,t}}$\par
$\kn{13a}{3641}$\par
$\kn{13a}{3660}$$^{\mathrm{d}}$\par
$\kn{13a}{3718}$\par
$\kn{13a}{3749}$$^{\mathrm{d}}$\par
$\kn{13a}{3794}$\par
$\kn{13a}{3796}$$^{\mathrm{d}}$\par
$\kn{13a}{3841}$$^{\mathrm{d}}$\par
$\kn{13a}{3865}$\par
$\kn{13a}{3880}$$^{\mathrm{d}}$\par
$\kn{13a}{3913}$\par
$\kn{13a}{3917}$\par
$\kn{13a}{3931}$\par
$\kn{13a}{3978}$$^{\mathrm{d}}$\par
$\kn{13a}{4017}$\par
$\kn{13a}{4039}$$^{\mathrm{d,t}}$\par
$\kn{13a}{4047}$$^{\mathrm{d}}$\par
$\kn{13a}{4055}$\par
$\kn{13a}{4112}$\par
$\kn{13a}{4127}$$^{\mathrm{d,t}}$\par
$\kn{13a}{4176}$\par
$\kn{13a}{4182}$\par
$\kn{13a}{4210}$\par
$\kn{13a}{4250}$\par
$\kn{13a}{4307}$$^{\mathrm{d}}$\par
$\kn{13a}{4359}$\par
$\kn{13a}{4578}$$^{\mathrm{d,t}}$\par
$\kn{13a}{4631}$$^{\mathrm{d}}$\par
$\kn{13a}{4695}$\par
$\kn{13a}{4706}$\par
$\kn{13a}{4745}$\par
$\kn{13n}{2}$$^{\mathrm{p}}$\par
$\kn{13n}{16}$$^{\mathrm{p}}$\par
$\kn{13n}{139}$$^{\mathrm{p}}$\par
$\kn{13n}{385}$\par
$\kn{13n}{404}$\par
$\kn{13n}{416}$\par
$\kn{13n}{439}$\par
$\kn{13n}{446}$$^{\mathrm{D}}$\par
$\kn{13n}{448}$\par
$\kn{13n}{472}$$^{\mathrm{s,t}}$\par
$\kn{13n}{554}$$^{\mathrm{s}}$\par
$\kn{13n}{560}$$^{\mathrm{p}}$\par
$\kn{13n}{582}$$^{\mathrm{D}}$\par
$\kn{13n}{628}$$^{\mathrm{p,t}}$\par
$\kn{13n}{638}$$^{\mathrm{p}}$\par
$\kn{13n}{650}$\par
$\kn{13n}{674}$$^{\mathrm{p}}$\par
$\kn{13n}{704}$$^{\mathrm{D}}$\par
$\kn{13n}{780}$$^{\mathrm{D}}$\par
$\kn{13n}{790}$$^{\mathrm{p}}$\par
$\kn{13n}{846}$$^{\mathrm{D}}$\par
$\kn{13n}{899}$$^{\mathrm{D}}$\par
$\kn{13n}{926}$$^{\mathrm{D}}$\par
$\kn{13n}{946}$$^{\mathrm{D}}$\par
$\kn{13n}{966}$$^{\mathrm{p}}$\par
$\kn{13n}{969}$$^{\mathrm{D}}$\par
$\kn{13n}{984}$\par
$\kn{13n}{990}$\par
$\kn{13n}{1000}$$^{\mathrm{D}}$\par
$\kn{13n}{1002}$$^{\mathrm{p}}$\par
$\kn{13n}{1042}$$^{\mathrm{D}}$\par
$\kn{13n}{1055}$\par
$\kn{13n}{1078}$$^{\mathrm{p}}$\par
$\kn{13n}{1134}$$^{\mathrm{p,t}}$\par
$\kn{13n}{1157}$\par
$\kn{13n}{1161}$$^{\mathrm{p}}$\par
$\kn{13n}{1183}$$^{\mathrm{p}}$\par
$\kn{13n}{1193}$\par
$\kn{13n}{1197}$\par
$\kn{13n}{1281}$$^{\mathrm{p,t}}$\par
$\kn{13n}{1307}$$^{\mathrm{p,t}}$\par
$\kn{13n}{1346}$$^{\mathrm{D}}$\par
$\kn{13n}{1377}$\par
$\kn{13n}{1379}$\par
$\kn{13n}{1469}$$^{\mathrm{D}}$\par
$\kn{13n}{1518}$\par
$\kn{13n}{1530}$\par
$\kn{13n}{1553}$$^{\mathrm{D}}$\par
$\kn{13n}{1576}$$^{\mathrm{p,t}}$\par
$\kn{13n}{1618}$$^{\mathrm{p,t}}$\par
$\kn{13n}{1621}$$^{\mathrm{D}}$\par
$\kn{13n}{1861}$$^{\mathrm{s,t}}$\par
$\kn{13n}{2010}$$^{\mathrm{s,t}}$\par
$\kn{13n}{2032}$$^{\mathrm{s,t}}$\par
$\kn{13n}{2067}$$^{\mathrm{p}}$\par
$\kn{13n}{2182}$$^{\mathrm{p,t}}$\par
$\kn{13n}{2206}$$^{\mathrm{s,t}}$\par
$\kn{13n}{2222}$$^{\mathrm{D}}$\par
$\kn{13n}{2276}$$^{\mathrm{p}}$\par
$\kn{13n}{2326}$\par
$\kn{13n}{2341}$\par
$\kn{13n}{2352}$$^{\mathrm{p,t}}$\par
$\kn{13n}{2366}$$^{\mathrm{p,t}}$\par
$\kn{13n}{2454}$\par
$\kn{13n}{2510}$$^{\mathrm{D}}$\par
$\kn{13n}{2574}$\par
$\kn{13n}{2626}$$^{\mathrm{D}}$\par
$\kn{13n}{2632}$$^{\mathrm{D}}$\par
$\kn{13n}{2654}$\par
$\kn{13n}{2670}$\par
$\kn{13n}{2675}$\par
$\kn{13n}{2683}$$^{\mathrm{D}}$\par
$\kn{13n}{2764}$\par
$\kn{13n}{2811}$\par
$\kn{13n}{2850}$\par
$\kn{13n}{2893}$$^{\mathrm{p}}$\par
$\kn{13n}{2931}$\par
$\kn{13n}{2932}$$^{\mathrm{p,t}}$\par
$\kn{13n}{2998}$\par
$\kn{13n}{3093}$$^{\mathrm{D}}$\par
$\kn{13n}{3183}$$^{\mathrm{s,t}}$\par
$\kn{13n}{3232}$$^{\mathrm{D}}$\par
$\kn{13n}{3234}$$^{\mathrm{s,t}}$\par
$\kn{13n}{3249}$$^{\mathrm{p,t}}$\par
$\kn{13n}{3252}$$^{\mathrm{D}}$\par
$\kn{13n}{3266}$\par
$\kn{13n}{3271}$\par
$\kn{13n}{3283}$$^{\mathrm{D}}$\par
$\kn{13n}{3309}$$^{\mathrm{D}}$\par
$\kn{13n}{3430}$\par
$\kn{13n}{3439}$$^{\mathrm{D}}$\par
$\kn{13n}{3458}$\par
$\kn{13n}{3512}$\par
$\kn{13n}{3530}$\par
$\kn{13n}{3547}$$^{\mathrm{p,t}}$\par
$\kn{13n}{3597}$$^{\mathrm{D}}$\par
$\kn{13n}{3599}$$^{\mathrm{D}}$\par
$\kn{13n}{3790}$\par
$\kn{13n}{3794}$$^{\mathrm{D}}$\par
$\kn{13n}{3844}$$^{\mathrm{s,t}}$\par
$\kn{13n}{3855}$$^{\mathrm{s,t}}$\par
$\kn{13n}{3884}$$^{\mathrm{D}}$\par
$\kn{13n}{4111}$$^{\mathrm{s}}$\par
$\kn{13n}{4171}$\par
$\kn{13n}{4199}$\par
$\kn{13n}{4292}$$^{\mathrm{D}}$\par
$\kn{13n}{4308}$$^{\mathrm{s}}$\par
$\kn{13n}{4369}$\par
$\kn{13n}{4396}$\par
$\kn{13n}{4507}$\par
$\kn{13n}{4748}$\par
$\kn{13n}{4764}$\par
$\kn{13n}{4798}$$^{\mathrm{D}}$\par
$\kn{13n}{4815}$\par
$\kn{13n}{4919}$$^{\mathrm{D}}$\par
$\kn{13n}{4947}$$^{\mathrm{s,t}}$\par
$\kn{13n}{4955}$$^{\mathrm{D}}$\par
$\kn{13n}{5035}$$^{\mathrm{D}}$\par
\end{multicols}
\noindent\begin{minipage}{\textwidth}
\hrule height\heavyrulewidth\vspace{3pt}
\raggedright\setlength{\parindent}{0pt}\setlength{\parskip}{0pt}
\strut Unmarked: $h_\lambda=2$, so $\gtop=\gfour=2$ (Corollary~\ref{cor:criteria}(b)).\strut\par
\strut $^{\mathrm{d}}$ $h_\lambda=1$, and the correction terms of the Goeritz manifold or of the star plumbing exclude a smooth M\"obius band (Proposition~\ref{prop:lower}(c)).\strut\par
\strut $^{\mathrm{p}}$, $^{\mathrm{s}}$ the same, with the correction terms from a plumbing tree or from a surgery description (Proposition~\ref{prop:newd}(a), (b)).\strut\par
\strut $^{\mathrm{D}}$ $h_\lambda=1$, and Donaldson's theorem, applied to a definite checkerboard surface, excludes a smooth M\"obius band (Proposition~\ref{prop:newd}(c)).\strut\par
\strut $^{\mathrm{t}}$ moreover $\gtop=1$ (Theorem~\ref{mthm:gap}).\strut\par
\end{minipage}\par}

\Needspace{8\baselineskip}
\subsection{Knots with \texorpdfstring{$\gfour=2$}{gamma4 = 2} and range \texorpdfstring{$[1,3]$}{[1,3]}}\label{app:newthree}
The following $\nGammaTwoOneThree$ knots, whose \KI{} range is $[1,3]$, have $\gfour=2$ (\S\ref{sec:values}).  The entries are ordered down the columns.

{\footnotesize
\setlength{\columnsep}{4pt}
\setlength{\multicolsep}{3pt}
\par\vspace{3pt}\hrule height\heavyrulewidth
\begin{multicols}{7}
\raggedright\setlength{\parindent}{0pt}\setlength{\parskip}{0pt}
$\kn{11a}{203}$$^{\mathrm{d}}$\par
$\kn{11a}{206}$$^{\mathrm{d}}$\par
$\kn{11a}{259}$$^{\mathrm{d}}$\par
$\kn{11a}{297}$\par
$\kn{11a}{305}$\par
$\kn{11a}{308}$$^{\mathrm{d}}$\par
$\kn{11a}{314}$\par
$\kn{11a}{318}$\par
$\kn{11n}{57}$$^{\mathrm{d}}$\par
$\kn{12a}{36}$$^{\mathrm{d}}$\par
$\kn{12a}{38}$$^{\mathrm{d}}$\par
$\kn{12a}{156}$$^{\mathrm{d}}$\par
$\kn{12a}{176}$$^{\mathrm{d}}$\par
$\kn{12a}{195}$$^{\mathrm{d}}$\par
$\kn{12a}{332}$\par
$\kn{12a}{493}$$^{\mathrm{d}}$\par
$\kn{12a}{587}$$^{\mathrm{d}}$\par
$\kn{12a}{590}$$^{\mathrm{d}}$\par
$\kn{12a}{647}$\par
$\kn{12a}{669}$$^{\mathrm{d}}$\par
$\kn{12a}{693}$$^{\mathrm{d}}$\par
$\kn{12a}{694}$$^{\mathrm{d}}$\par
$\kn{12a}{714}$$^{\mathrm{d}}$\par
$\kn{12a}{726}$$^{\mathrm{d}}$\par
$\kn{12a}{742}$\par
$\kn{12a}{801}$$^{\mathrm{d}}$\par
$\kn{12a}{802}$$^{\mathrm{d}}$\par
$\kn{12a}{827}$$^{\mathrm{d}}$\par
$\kn{12a}{829}$$^{\mathrm{d}}$\par
$\kn{12a}{832}$$^{\mathrm{d}}$\par
$\kn{12a}{840}$$^{\mathrm{d}}$\par
$\kn{12a}{853}$$^{\mathrm{d}}$\par
$\kn{12a}{855}$$^{\mathrm{d}}$\par
$\kn{12a}{859}$$^{\mathrm{d}}$\par
$\kn{12a}{923}$\par
$\kn{12a}{937}$$^{\mathrm{d}}$\par
$\kn{12a}{946}$$^{\mathrm{d}}$\par
$\kn{12a}{984}$$^{\mathrm{d}}$\par
$\kn{12a}{1012}$$^{\mathrm{d}}$\par
$\kn{12a}{1142}$\par
$\kn{12a}{1171}$$^{\mathrm{d}}$\par
$\kn{12a}{1179}$$^{\mathrm{d}}$\par
$\kn{12a}{1260}$\par
$\kn{12n}{163}$$^{\mathrm{d}}$\par
$\kn{12n}{172}$$^{\mathrm{d}}$\par
$\kn{12n}{270}$\par
$\kn{12n}{293}$$^{\mathrm{d}}$\par
$\kn{12n}{494}$\par
$\kn{12n}{505}$\par
$\kn{12n}{518}$\par
$\kn{12n}{546}$\par
$\kn{12n}{571}$$^{\mathrm{d}}$\par
$\kn{12n}{598}$\par
$\kn{12n}{602}$\par
$\kn{12n}{611}$\par
$\kn{12n}{654}$\par
$\kn{12n}{701}$\par
$\kn{12n}{756}$\par
$\kn{12n}{806}$\par
$\kn{12n}{844}$\par
$\kn{12n}{847}$\par
$\kn{12n}{850}$\par
$\kn{13a}{125}$$^{\mathrm{d}}$\par
$\kn{13a}{132}$\par
$\kn{13a}{133}$$^{\mathrm{d}}$\par
$\kn{13a}{173}$$^{\mathrm{d}}$\par
$\kn{13a}{191}$\par
$\kn{13a}{193}$$^{\mathrm{d}}$\par
$\kn{13a}{195}$$^{\mathrm{d}}$\par
$\kn{13a}{254}$$^{\mathrm{d}}$\par
$\kn{13a}{279}$$^{\mathrm{d}}$\par
$\kn{13a}{286}$$^{\mathrm{d}}$\par
$\kn{13a}{382}$$^{\mathrm{d}}$\par
$\kn{13a}{416}$$^{\mathrm{d}}$\par
$\kn{13a}{435}$$^{\mathrm{d}}$\par
$\kn{13a}{448}$$^{\mathrm{d}}$\par
$\kn{13a}{469}$$^{\mathrm{d}}$\par
$\kn{13a}{501}$$^{\mathrm{d}}$\par
$\kn{13a}{561}$\par
$\kn{13a}{616}$$^{\mathrm{d}}$\par
$\kn{13a}{627}$$^{\mathrm{d}}$\par
$\kn{13a}{648}$$^{\mathrm{d}}$\par
$\kn{13a}{652}$$^{\mathrm{d}}$\par
$\kn{13a}{665}$$^{\mathrm{d}}$\par
$\kn{13a}{667}$$^{\mathrm{d}}$\par
$\kn{13a}{681}$$^{\mathrm{d}}$\par
$\kn{13a}{705}$\par
$\kn{13a}{720}$\par
$\kn{13a}{855}$\par
$\kn{13a}{896}$\par
$\kn{13a}{905}$$^{\mathrm{d}}$\par
$\kn{13a}{1043}$$^{\mathrm{d}}$\par
$\kn{13a}{1133}$\par
$\kn{13a}{1157}$$^{\mathrm{d}}$\par
$\kn{13a}{1196}$\par
$\kn{13a}{1243}$$^{\mathrm{d}}$\par
$\kn{13a}{1335}$\par
$\kn{13a}{1353}$\par
$\kn{13a}{1462}$$^{\mathrm{d}}$\par
$\kn{13a}{1467}$$^{\mathrm{d}}$\par
$\kn{13a}{1468}$\par
$\kn{13a}{1482}$\par
$\kn{13a}{1492}$$^{\mathrm{d}}$\par
$\kn{13a}{1520}$$^{\mathrm{d}}$\par
$\kn{13a}{1617}$\par
$\kn{13a}{1645}$$^{\mathrm{d}}$\par
$\kn{13a}{1725}$\par
$\kn{13a}{1740}$$^{\mathrm{d}}$\par
$\kn{13a}{1929}$$^{\mathrm{d}}$\par
$\kn{13a}{1969}$\par
$\kn{13a}{1986}$\par
$\kn{13a}{1990}$$^{\mathrm{d}}$\par
$\kn{13a}{2020}$$^{\mathrm{d}}$\par
$\kn{13a}{2038}$\par
$\kn{13a}{2097}$$^{\mathrm{d}}$\par
$\kn{13a}{2124}$$^{\mathrm{d}}$\par
$\kn{13a}{2129}$$^{\mathrm{d}}$\par
$\kn{13a}{2154}$$^{\mathrm{d}}$\par
$\kn{13a}{2157}$$^{\mathrm{d}}$\par
$\kn{13a}{2173}$$^{\mathrm{d}}$\par
$\kn{13a}{2175}$$^{\mathrm{d}}$\par
$\kn{13a}{2232}$\par
$\kn{13a}{2233}$\par
$\kn{13a}{2271}$\par
$\kn{13a}{2319}$$^{\mathrm{d}}$\par
$\kn{13a}{2362}$$^{\mathrm{d}}$\par
$\kn{13a}{2445}$$^{\mathrm{d}}$\par
$\kn{13a}{2507}$$^{\mathrm{d}}$\par
$\kn{13a}{2548}$$^{\mathrm{d}}$\par
$\kn{13a}{2564}$$^{\mathrm{d}}$\par
$\kn{13a}{2620}$$^{\mathrm{d}}$\par
$\kn{13a}{2642}$\par
$\kn{13a}{2647}$$^{\mathrm{d}}$\par
$\kn{13a}{2667}$\par
$\kn{13a}{2719}$$^{\mathrm{d}}$\par
$\kn{13a}{2721}$$^{\mathrm{d}}$\par
$\kn{13a}{2810}$\par
$\kn{13a}{2824}$$^{\mathrm{d}}$\par
$\kn{13a}{2830}$$^{\mathrm{d}}$\par
$\kn{13a}{2833}$$^{\mathrm{d}}$\par
$\kn{13a}{2835}$$^{\mathrm{d}}$\par
$\kn{13a}{2844}$\par
$\kn{13a}{2879}$\par
$\kn{13a}{2884}$$^{\mathrm{d}}$\par
$\kn{13a}{2944}$$^{\mathrm{d}}$\par
$\kn{13a}{2945}$$^{\mathrm{d}}$\par
$\kn{13a}{3018}$$^{\mathrm{d}}$\par
$\kn{13a}{3110}$$^{\mathrm{d}}$\par
$\kn{13a}{3122}$\par
$\kn{13a}{3128}$\par
$\kn{13a}{3139}$\par
$\kn{13a}{3158}$\par
$\kn{13a}{3164}$$^{\mathrm{d}}$\par
$\kn{13a}{3206}$$^{\mathrm{d}}$\par
$\kn{13a}{3234}$$^{\mathrm{d}}$\par
$\kn{13a}{3250}$$^{\mathrm{d}}$\par
$\kn{13a}{3262}$$^{\mathrm{d}}$\par
$\kn{13a}{3282}$$^{\mathrm{d}}$\par
$\kn{13a}{3309}$$^{\mathrm{d}}$\par
$\kn{13a}{3315}$$^{\mathrm{d}}$\par
$\kn{13a}{3329}$$^{\mathrm{d}}$\par
$\kn{13a}{3337}$$^{\mathrm{d}}$\par
$\kn{13a}{3340}$$^{\mathrm{d}}$\par
$\kn{13a}{3358}$$^{\mathrm{d}}$\par
$\kn{13a}{3377}$$^{\mathrm{d}}$\par
$\kn{13a}{3380}$$^{\mathrm{d}}$\par
$\kn{13a}{3387}$$^{\mathrm{d}}$\par
$\kn{13a}{3459}$\par
$\kn{13a}{3461}$$^{\mathrm{d}}$\par
$\kn{13a}{3524}$$^{\mathrm{d}}$\par
$\kn{13a}{3530}$\par
$\kn{13a}{3554}$$^{\mathrm{d}}$\par
$\kn{13a}{3558}$$^{\mathrm{d}}$\par
$\kn{13a}{3589}$$^{\mathrm{d}}$\par
$\kn{13a}{3642}$\par
$\kn{13a}{3670}$$^{\mathrm{d}}$\par
$\kn{13a}{3687}$\par
$\kn{13a}{3693}$\par
$\kn{13a}{3706}$$^{\mathrm{d}}$\par
$\kn{13a}{3720}$$^{\mathrm{d}}$\par
$\kn{13a}{3800}$$^{\mathrm{d}}$\par
$\kn{13a}{3852}$$^{\mathrm{d}}$\par
$\kn{13a}{3858}$\par
$\kn{13a}{3892}$$^{\mathrm{d}}$\par
$\kn{13a}{3901}$$^{\mathrm{d}}$\par
$\kn{13a}{3950}$$^{\mathrm{d}}$\par
$\kn{13a}{3956}$\par
$\kn{13a}{3961}$$^{\mathrm{d}}$\par
$\kn{13a}{3990}$$^{\mathrm{d}}$\par
$\kn{13a}{4030}$\par
$\kn{13a}{4034}$\par
$\kn{13a}{4041}$$^{\mathrm{d}}$\par
$\kn{13a}{4053}$$^{\mathrm{d}}$\par
$\kn{13a}{4102}$$^{\mathrm{d}}$\par
$\kn{13a}{4116}$\par
$\kn{13a}{4129}$\par
$\kn{13a}{4163}$$^{\mathrm{d}}$\par
$\kn{13a}{4190}$\par
$\kn{13a}{4201}$\par
$\kn{13a}{4214}$$^{\mathrm{d}}$\par
$\kn{13a}{4240}$$^{\mathrm{d}}$\par
$\kn{13a}{4258}$$^{\mathrm{d}}$\par
$\kn{13a}{4306}$\par
$\kn{13a}{4309}$\par
$\kn{13a}{4313}$$^{\mathrm{d}}$\par
$\kn{13a}{4329}$\par
$\kn{13a}{4401}$$^{\mathrm{d}}$\par
$\kn{13a}{4427}$\par
$\kn{13a}{4488}$$^{\mathrm{d}}$\par
$\kn{13a}{4504}$$^{\mathrm{d}}$\par
$\kn{13a}{4509}$$^{\mathrm{d}}$\par
$\kn{13a}{4532}$$^{\mathrm{d}}$\par
$\kn{13a}{4538}$$^{\mathrm{d}}$\par
$\kn{13a}{4545}$$^{\mathrm{d}}$\par
$\kn{13a}{4553}$\par
$\kn{13a}{4558}$$^{\mathrm{d}}$\par
$\kn{13a}{4576}$\par
$\kn{13a}{4595}$$^{\mathrm{d}}$\par
$\kn{13a}{4607}$$^{\mathrm{d}}$\par
$\kn{13a}{4616}$\par
$\kn{13a}{4643}$$^{\mathrm{d}}$\par
$\kn{13a}{4658}$$^{\mathrm{d}}$\par
$\kn{13a}{4678}$$^{\mathrm{d}}$\par
$\kn{13a}{4698}$$^{\mathrm{d}}$\par
$\kn{13a}{4739}$\par
$\kn{13a}{4764}$\par
$\kn{13a}{4768}$\par
$\kn{13a}{4800}$\par
$\kn{13a}{4828}$$^{\mathrm{d}}$\par
$\kn{13a}{4844}$\par
$\kn{13a}{4856}$$^{\mathrm{d}}$\par
$\kn{13a}{4857}$\par
$\kn{13a}{4862}$$^{\mathrm{d}}$\par
$\kn{13a}{4876}$\par
$\kn{13n}{12}$$^{\mathrm{p}}$\par
$\kn{13n}{19}$$^{\mathrm{p}}$\par
$\kn{13n}{85}$\par
$\kn{13n}{152}$\par
$\kn{13n}{166}$$^{\mathrm{p}}$\par
$\kn{13n}{178}$$^{\mathrm{p}}$\par
$\kn{13n}{199}$$^{\mathrm{p}}$\par
$\kn{13n}{218}$$^{\mathrm{p}}$\par
$\kn{13n}{277}$$^{\mathrm{p}}$\par
$\kn{13n}{383}$$^{\mathrm{p}}$\par
$\kn{13n}{401}$$^{\mathrm{D}}$\par
$\kn{13n}{460}$$^{\mathrm{D}}$\par
$\kn{13n}{479}$$^{\mathrm{D}}$\par
$\kn{13n}{488}$$^{\mathrm{D}}$\par
$\kn{13n}{537}$$^{\mathrm{D}}$\par
$\kn{13n}{573}$$^{\mathrm{p}}$\par
$\kn{13n}{594}$$^{\mathrm{p,t}}$\par
$\kn{13n}{652}$$^{\mathrm{D}}$\par
$\kn{13n}{694}$$^{\mathrm{D}}$\par
$\kn{13n}{697}$$^{\mathrm{p}}$\par
$\kn{13n}{828}$\par
$\kn{13n}{831}$\par
$\kn{13n}{874}$\par
$\kn{13n}{949}$\par
$\kn{13n}{972}$$^{\mathrm{D}}$\par
$\kn{13n}{1003}$$^{\mathrm{p,t}}$\par
$\kn{13n}{1079}$\par
$\kn{13n}{1148}$$^{\mathrm{p}}$\par
$\kn{13n}{1164}$$^{\mathrm{p,t}}$\par
$\kn{13n}{1165}$$^{\mathrm{p}}$\par
$\kn{13n}{1234}$$^{\mathrm{D}}$\par
$\kn{13n}{1297}$$^{\mathrm{p}}$\par
$\kn{13n}{1380}$$^{\mathrm{p}}$\par
$\kn{13n}{1385}$$^{\mathrm{p}}$\par
$\kn{13n}{1409}$\par
$\kn{13n}{1411}$$^{\mathrm{p,t}}$\par
$\kn{13n}{1437}$$^{\mathrm{p,t}}$\par
$\kn{13n}{1460}$\par
$\kn{13n}{1489}$$^{\mathrm{p}}$\par
$\kn{13n}{1506}$$^{\mathrm{p}}$\par
$\kn{13n}{1524}$\par
$\kn{13n}{1557}$\par
$\kn{13n}{1747}$$^{\mathrm{p}}$\par
$\kn{13n}{1923}$\par
$\kn{13n}{1977}$$^{\mathrm{p}}$\par
$\kn{13n}{1987}$$^{\mathrm{p}}$\par
$\kn{13n}{2009}$\par
$\kn{13n}{2013}$$^{\mathrm{p}}$\par
$\kn{13n}{2155}$\par
$\kn{13n}{2271}$$^{\mathrm{p,t}}$\par
$\kn{13n}{2312}$\par
$\kn{13n}{2386}$$^{\mathrm{p,t}}$\par
$\kn{13n}{2391}$$^{\mathrm{p}}$\par
$\kn{13n}{2402}$\par
$\kn{13n}{2455}$\par
$\kn{13n}{2465}$$^{\mathrm{D}}$\par
$\kn{13n}{2469}$\par
$\kn{13n}{2486}$$^{\mathrm{D}}$\par
$\kn{13n}{2516}$\par
$\kn{13n}{2550}$\par
$\kn{13n}{2565}$$^{\mathrm{p}}$\par
$\kn{13n}{2612}$\par
$\kn{13n}{2732}$$^{\mathrm{D}}$\par
$\kn{13n}{2860}$\par
$\kn{13n}{2862}$\par
$\kn{13n}{2909}$\par
$\kn{13n}{2925}$\par
$\kn{13n}{2994}$\par
$\kn{13n}{3004}$\par
$\kn{13n}{3020}$\par
$\kn{13n}{3089}$$^{\mathrm{D}}$\par
$\kn{13n}{3133}$$^{\mathrm{D}}$\par
$\kn{13n}{3166}$\par
$\kn{13n}{3224}$$^{\mathrm{D}}$\par
$\kn{13n}{3263}$\par
$\kn{13n}{3303}$$^{\mathrm{D}}$\par
$\kn{13n}{3330}$$^{\mathrm{D}}$\par
$\kn{13n}{3388}$$^{\mathrm{D}}$\par
$\kn{13n}{3507}$\par
$\kn{13n}{3529}$\par
$\kn{13n}{3587}$$^{\mathrm{D}}$\par
$\kn{13n}{3603}$$^{\mathrm{D}}$\par
$\kn{13n}{3613}$$^{\mathrm{D}}$\par
$\kn{13n}{3616}$$^{\mathrm{D}}$\par
$\kn{13n}{3639}$\par
$\kn{13n}{3659}$\par
$\kn{13n}{3718}$\par
$\kn{13n}{3720}$$^{\mathrm{D}}$\par
$\kn{13n}{3732}$$^{\mathrm{D}}$\par
$\kn{13n}{3891}$\par
$\kn{13n}{3918}$\par
$\kn{13n}{3929}$\par
$\kn{13n}{4013}$$^{\mathrm{D}}$\par
$\kn{13n}{4173}$\par
$\kn{13n}{4225}$\par
$\kn{13n}{4307}$$^{\mathrm{D}}$\par
$\kn{13n}{4381}$\par
$\kn{13n}{4497}$\par
$\kn{13n}{4672}$\par
$\kn{13n}{4698}$$^{\mathrm{D}}$\par
$\kn{13n}{4817}$$^{\mathrm{D}}$\par
$\kn{13n}{4823}$$^{\mathrm{D}}$\par
$\kn{13n}{4964}$\par
$\kn{13n}{4971}$\par
$\kn{13n}{5027}$\par
$\kn{13n}{5060}$\par
$\kn{13n}{5062}$\par
$\kn{13n}{5107}$\par
\end{multicols}
\noindent\begin{minipage}{\textwidth}
\hrule height\heavyrulewidth\vspace{3pt}
\raggedright\setlength{\parindent}{0pt}\setlength{\parskip}{0pt}
\strut $^{\mathrm{d}}$ $h_\lambda=1$, and the correction terms of the Goeritz manifold or of the star plumbing exclude a smooth M\"obius band (Proposition~\ref{prop:lower}(c)).\strut\par
\strut $^{\mathrm{p}}$ the same, with the correction terms from a plumbing tree (Proposition~\ref{prop:newd}(a)).\strut\par
\strut $^{\mathrm{D}}$ $h_\lambda=1$, and Donaldson's theorem, applied to a definite checkerboard surface, excludes a smooth M\"obius band (Proposition~\ref{prop:newd}(c)).\strut\par
\strut $^{\mathrm{t}}$ moreover $\gtop=1$ (Theorem~\ref{mthm:gap}).\strut\par
\strut Unmarked knots have $h_\lambda=2$, so $\gtop=\gfour=2$; every knot listed bounds a punctured Klein bottle made of band moves (Proposition~\ref{prop:upper}(b)).\strut\par
\end{minipage}\par}

\Needspace{14\baselineskip}
\subsection{Knots with \texorpdfstring{$h_\lambda=3$}{h = 3}}
Table~\ref{tab:hthree} lists the $\nHthree$ knots with at most $13$ crossings and $h_\lambda=3$, with the invariants of
Corollary~\ref{cor:criteria}(c) that give $h_\lambda=3$.  All of them have $\gfour=\gtop=3$: the $\nGammaThreeAll$ knots with $\gfour=3$
of \S\ref{sec:values}, $\kn{8}{18}$, whose value $3$ \KI{} lists, and $\kn{9}{49}$ and $\kn{11n}{133}$
(\S\ref{sec:corrections}).

{\small\setlength{\tabcolsep}{4pt}
\begin{longtable}{lrlrrlcccc}
\caption{The $49$ knots with at most $13$ crossings and $h_\lambda=3$, with the invariants of Corollary~\ref{cor:criteria}(c): $s=\sigma+4\Arf\bmod 8$, the primes $p$ with $a_p=2$ and their signs $\varepsilon_p$, and the conditions (i)--(iii) that hold.  The column $\ws$ gives the rational Witt span \cite[Thm.~4.6]{LeeTheory}.  The last two columns give the value or range of $\gfour$ in \KI{} 2026.11.2 and the value of $\gfour$; it was known for $\kn{8}{18}$ (\KI{}) and $\kn{9}{49}$ (the revised version \cite{JKv3} of \cite{JK}), and is determined in this paper for the others.}\label{tab:hthree}\\
\toprule
Knot & $\det$ & $H_K$ & $\sigma$ & $s$ & $p{:}\,\varepsilon_p$ & criterion & $\ws$ & \KI{} & value \\
\midrule
\endfirsthead
\multicolumn{10}{@{}l@{}}{\emph{Knots with $h_\lambda=3$, continued from the previous page.}}\\[2pt]
\toprule
Knot & $\det$ & $H_K$ & $\sigma$ & $s$ & $p{:}\,\varepsilon_p$ & criterion & $\ws$ & \KI{} & value \\
\midrule
\endhead
\bottomrule
\multicolumn{10}{@{}r@{}}{\emph{continued on the next page}}\\
\endfoot
\bottomrule
\endlastfoot
$\kn{8}{18}$ & $45$ & $\Z/3 \oplus \Z/3 \oplus \Z/5$ & $0$ & $4$ & $3{:}\,-$ & (iii) & $2$ & $3$ & $3$ \\
$\kn{9}{49}$ & $25$ & $\Z/5 \oplus \Z/5$ & $-4$ & $4$ & $5{:}\,-$ & (i) & $3$ & $2$ & $3$ \\
$\kn{11n}{133}$ & $25$ & $\Z/5 \oplus \Z/5$ & $-4$ & $4$ & $5{:}\,-$ & (i) & $3$ & $2$ & $3$ \\
$\kn{12a}{561}$ & $225$ & $\Z/9 \oplus \Z/5 \oplus \Z/5$ & $-4$ & $4$ & $5{:}\,-$ & (i) & $3$ & $[2,3]$ & $3$ \\
$\kn{12a}{569}$ & $153$ & $\Z/3 \oplus \Z/3 \oplus \Z/17$ & $4$ & $4$ & $3{:}\,-$ & (iii) & $3$ & $[2,3]$ & $3$ \\
$\kn{12a}{664}$ & $169$ & $\Z/13 \oplus \Z/13$ & $-4$ & $4$ & $13{:}\,-$ & (i) & $3$ & $[2,3]$ & $3$ \\
$\kn{12a}{780}$ & $225$ & $\Z/9 \oplus \Z/5 \oplus \Z/5$ & $-4$ & $4$ & $5{:}\,-$ & (i) & $3$ & $[2,3]$ & $3$ \\
$\kn{12a}{906}$ & $245$ & $\Z/5 \oplus \Z/7 \oplus \Z/7$ & $0$ & $4$ & $7{:}\,-$ & (iii) & $2$ & $[2,3]$ & $3$ \\
$\kn{12a}{907}$ & $225$ & $\Z/9 \oplus \Z/5 \oplus \Z/5$ & $-4$ & $4$ & $5{:}\,-$ & (i) & $3$ & $[2,3]$ & $3$ \\
$\kn{12a}{1123}$ & $261$ & $\Z/3 \oplus \Z/3 \oplus \Z/29$ & $0$ & $4$ & $3{:}\,-$ & (iii) & $2$ & $[2,3]$ & $3$ \\
$\kn{12n}{276}$ & $25$ & $\Z/5 \oplus \Z/5$ & $4$ & $4$ & $5{:}\,-$ & (i) & $3$ & $[2,3]$ & $3$ \\
$\kn{12n}{333}$ & $45$ & $\Z/3 \oplus \Z/3 \oplus \Z/5$ & $0$ & $4$ & $3{:}\,-$ & (iii) & $2$ & $[2,3]$ & $3$ \\
$\kn{12n}{626}$ & $153$ & $\Z/3 \oplus \Z/3 \oplus \Z/17$ & $4$ & $4$ & $3{:}\,-$ & (iii) & $3$ & $[2,3]$ & $3$ \\
$\kn{12n}{888}$ & $45$ & $\Z/3 \oplus \Z/3 \oplus \Z/5$ & $-8$ & $4$ & $3{:}\,-$ & (iii) & $2$ & $[2,3]$ & $3$ \\
$\kn{13a}{1398}$ & $225$ & $\Z/9 \oplus \Z/5 \oplus \Z/5$ & $-4$ & $4$ & $5{:}\,-$ & (i) & $3$ & $[2,3]$ & $3$ \\
$\kn{13a}{1557}$ & $225$ & $\Z/9 \oplus \Z/5 \oplus \Z/5$ & $4$ & $4$ & $5{:}\,-$ & (i) & $3$ & $[2,3]$ & $3$ \\
$\kn{13a}{1662}$ & $405$ & $\Z/3 \oplus \Z/27 \oplus \Z/5$ & $0$ & $4$ & $3{:}\,-$ & (iii) & $2$ & $[2,3]$ & $3$ \\
$\kn{13a}{1826}$ & $225$ & $\Z/9 \oplus \Z/5 \oplus \Z/5$ & $-4$ & $4$ & $5{:}\,-$ & (i) & $3$ & $[2,3]$ & $3$ \\
$\kn{13a}{1912}$ & $275$ & $\Z/5 \oplus \Z/5 \oplus \Z/11$ & $-6$ & $6$ & $5{:}\,-$ & (i) & $3$ & $[1,3]$ & $3$ \\
$\kn{13a}{1974}$ & $153$ & $\Z/3 \oplus \Z/3 \oplus \Z/17$ & $4$ & $4$ & $3{:}\,-$ & (iii) & $3$ & $[2,3]$ & $3$ \\
$\kn{13a}{1992}$ & $405$ & $\Z/3 \oplus \Z/27 \oplus \Z/5$ & $0$ & $4$ & $3{:}\,-$ & (iii) & $2$ & $[2,3]$ & $3$ \\
$\kn{13a}{2079}$ & $275$ & $\Z/5 \oplus \Z/5 \oplus \Z/11$ & $2$ & $6$ & $5{:}\,-$ & (i) & $3$ & $[1,3]$ & $3$ \\
$\kn{13a}{2181}$ & $275$ & $\Z/5 \oplus \Z/5 \oplus \Z/11$ & $-2$ & $2$ & $5{:}\,-$ & (i) & $3$ & $[1,3]$ & $3$ \\
$\kn{13a}{2890}$ & $289$ & $\Z/17 \oplus \Z/17$ & $-4$ & $4$ & $17{:}\,-$ & (i) & $3$ & $[2,4]$ & $3$ \\
$\kn{13a}{3034}$ & $405$ & $\Z/3 \oplus \Z/27 \oplus \Z/5$ & $0$ & $4$ & $3{:}\,-$ & (iii) & $2$ & $[2,3]$ & $3$ \\
$\kn{13a}{3118}$ & $475$ & $\Z/5 \oplus \Z/5 \oplus \Z/19$ & $-2$ & $2$ & $5{:}\,-$ & (i) & $3$ & $[1,3]$ & $3$ \\
$\kn{13a}{3121}$ & $153$ & $\Z/3 \oplus \Z/3 \oplus \Z/17$ & $-4$ & $4$ & $3{:}\,-$ & (iii) & $3$ & $[2,3]$ & $3$ \\
$\kn{13a}{3585}$ & $477$ & $\Z/3 \oplus \Z/3 \oplus \Z/53$ & $0$ & $4$ & $3{:}\,-$ & (iii) & $2$ & $[2,3]$ & $3$ \\
$\kn{13a}{3747}$ & $245$ & $\Z/5 \oplus \Z/7 \oplus \Z/7$ & $0$ & $4$ & $7{:}\,-$ & (iii) & $2$ & $[2,3]$ & $3$ \\
$\kn{13a}{3936}$ & $245$ & $\Z/5 \oplus \Z/7 \oplus \Z/7$ & $0$ & $4$ & $7{:}\,-$ & (iii) & $2$ & $[2,3]$ & $3$ \\
$\kn{13a}{4566}$ & $261$ & $\Z/3 \oplus \Z/3 \oplus \Z/29$ & $-8$ & $4$ & $3{:}\,-$ & (iii) & $2$ & $[2,3]$ & $3$ \\
$\kn{13a}{4637}$ & $405$ & $\Z/3 \oplus \Z/27 \oplus \Z/5$ & $0$ & $4$ & $3{:}\,-$ & (iii) & $2$ & $[2,3]$ & $3$ \\
$\kn{13a}{4674}$ & $405$ & $\Z/3 \oplus \Z/27 \oplus \Z/5$ & $0$ & $4$ & $3{:}\,-$ & (iii) & $2$ & $[2,3]$ & $3$ \\
$\kn{13n}{115}$ & $261$ & $\Z/3 \oplus \Z/3 \oplus \Z/29$ & $0$ & $4$ & $3{:}\,-$ & (iii) & $2$ & $[2,3]$ & $3$ \\
$\kn{13n}{210}$ & $45$ & $\Z/3 \oplus \Z/3 \oplus \Z/5$ & $0$ & $4$ & $3{:}\,-$ & (iii) & $2$ & $[2,3]$ & $3$ \\
$\kn{13n}{2305}$ & $25$ & $\Z/5 \oplus \Z/5$ & $4$ & $4$ & $5{:}\,-$ & (i) & $3$ & $[2,3]$ & $3$ \\
$\kn{13n}{2526}$ & $25$ & $\Z/5 \oplus \Z/5$ & $-4$ & $4$ & $5{:}\,-$ & (i) & $3$ & $[2,3]$ & $3$ \\
$\kn{13n}{3005}$ & $153$ & $\Z/3 \oplus \Z/3 \oplus \Z/17$ & $-4$ & $4$ & $3{:}\,-$ & (iii) & $3$ & $[2,3]$ & $3$ \\
$\kn{13n}{3253}$ & $169$ & $\Z/13 \oplus \Z/13$ & $-4$ & $4$ & $13{:}\,-$ & (i) & $3$ & $[2,3]$ & $3$ \\
$\kn{13n}{3446}$ & $169$ & $\Z/13 \oplus \Z/13$ & $-4$ & $4$ & $13{:}\,-$ & (i) & $3$ & $[2,3]$ & $3$ \\
$\kn{13n}{3489}$ & $25$ & $\Z/5 \oplus \Z/5$ & $-4$ & $4$ & $5{:}\,-$ & (i) & $3$ & $[2,3]$ & $3$ \\
$\kn{13n}{3546}$ & $25$ & $\Z/5 \oplus \Z/5$ & $-4$ & $4$ & $5{:}\,-$ & (i) & $3$ & $[2,3]$ & $3$ \\
$\kn{13n}{3570}$ & $169$ & $\Z/13 \oplus \Z/13$ & $-4$ & $4$ & $13{:}\,-$ & (i) & $3$ & $[2,3]$ & $3$ \\
$\kn{13n}{4696}$ & $25$ & $\Z/5 \oplus \Z/5$ & $4$ & $4$ & $5{:}\,-$ & (i) & $3$ & $[2,3]$ & $3$ \\
$\kn{13n}{4892}$ & $245$ & $\Z/5 \oplus \Z/7 \oplus \Z/7$ & $0$ & $4$ & $7{:}\,-$ & (iii) & $2$ & $[2,3]$ & $3$ \\
$\kn{13n}{4903}$ & $169$ & $\Z/13 \oplus \Z/13$ & $-4$ & $4$ & $13{:}\,-$ & (i) & $3$ & $[2,3]$ & $3$ \\
$\kn{13n}{4927}$ & $153$ & $\Z/3 \oplus \Z/3 \oplus \Z/17$ & $4$ & $4$ & $3{:}\,-$ & (iii) & $3$ & $[2,3]$ & $3$ \\
$\kn{13n}{5008}$ & $225$ & $\Z/9 \oplus \Z/5 \oplus \Z/5$ & $-4$ & $4$ & $5{:}\,-$ & (i) & $3$ & $[2,3]$ & $3$ \\
$\kn{13n}{5076}$ & $225$ & $\Z/3 \oplus \Z/3 \oplus \Z/5 \oplus \Z/5$ & $-4$ & $4$ & $5{:}\,-$ & (i) & $3$ & $[2,3]$ & $3$ \\
\end{longtable}}

\Needspace{8\baselineskip}
\subsection{Knots with fourteen crossings and \texorpdfstring{$h_\lambda=3$}{h = 3}}
The following $\nLSHthree$ knots of the table of \cite{LS} have $h_\lambda=3$ (Theorem~\ref{mthm:fourteen}).
The entries are ordered down the columns.

{\footnotesize
\setlength{\columnsep}{4pt}
\setlength{\multicolsep}{3pt}
\par\vspace{3pt}\hrule height\heavyrulewidth
\begin{multicols}{6}
\raggedright\setlength{\parindent}{0pt}\setlength{\parskip}{0pt}
\textup{K14a1463}\par
\textup{K14a1978}\par
\textup{K14a1987}\par
\textup{K14a2632}\par
\textup{K14a2920}\par
\textup{K14a2951}\par
\textup{K14a3382}\par
\textup{K14a4045}\par
\textup{K14a4631}\par
\textup{K14a4786}\par
\textup{K14a4874}\par
\textup{K14a5000}\par
\textup{K14a5638}\par
\textup{K14a5655}\par
\textup{K14a5668}\par
\textup{K14a6073}\par
\textup{K14a6401}\par
\textup{K14a6411}\par
\textup{K14a6683}\par
\textup{K14a6685}\par
\textup{K14a6697}\par
\textup{K14a6989}\par
\textup{K14a7183}\par
\textup{K14a7425}\par
\textup{K14a7498}\par
\textup{K14a7596}\par
\textup{K14a7914}\par
\textup{K14a8117}\par
\textup{K14a8654}\par
\textup{K14a8890}\par
\textup{K14a8910}\par
\textup{K14a8920}\par
\textup{K14a8948}\par
\textup{K14a9517}\par
\textup{K14a9533}\par
\textup{K14a9658}\par
\textup{K14a9950}\par
\textup{K14a10131}\par
\textup{K14a10886}\par
\textup{K14a10952}\par
\textup{K14a12669}\par
\textup{K14a12693}\par
\textup{K14a12759}\par
\textup{K14a12951}\par
\textup{K14a13528}\par
\textup{K14a14056}\par
\textup{K14a14258}\par
\textup{K14a14294}\par
\textup{K14a14533}\par
\textup{K14a14615}\par
\textup{K14a14959}\par
\textup{K14a15071}\par
\textup{K14a15141}\par
\textup{K14a15258}\par
\textup{K14a15290}\par
\textup{K14a15312}\par
\textup{K14a15447}\par
\textup{K14a16427}\par
\textup{K14a16801}\par
\textup{K14a17109}\par
\textup{K14a17382}\par
\textup{K14a17416}\par
\textup{K14a17442}\par
\textup{K14a17489}\par
\textup{K14a17761}\par
\textup{K14a17904}\par
\textup{K14a17926}\par
\textup{K14a17967}\par
\textup{K14a18153}\par
\textup{K14a18178}\par
\textup{K14a18353}\par
\textup{K14a18363}\par
\textup{K14a18545}\par
\textup{K14a18636}\par
\textup{K14a18853}\par
\textup{K14a18868}\par
\textup{K14a18975}\par
\textup{K14a18992}\par
\textup{K14a19052}\par
\textup{K14a19245}\par
\textup{K14a19270}\par
\textup{K14a19294}\par
\textup{K14a19476}\par
\textup{K14n568}\par
\textup{K14n1067}\par
\textup{K14n1371}\par
\textup{K14n1586}\par
\textup{K14n1662}\par
\textup{K14n1994}\par
\textup{K14n2008}\par
\textup{K14n2273}\par
\textup{K14n2283}\par
\textup{K14n2532}\par
\textup{K14n2710}\par
\textup{K14n3165}\par
\textup{K14n3243}\par
\textup{K14n3418}\par
\textup{K14n3491}\par
\textup{K14n3748}\par
\textup{K14n3823}\par
\textup{K14n3869}\par
\textup{K14n3872}\par
\textup{K14n3888}\par
\textup{K14n4007}\par
\textup{K14n4075}\par
\textup{K14n4398}\par
\textup{K14n4472}\par
\textup{K14n4514}\par
\textup{K14n4681}\par
\textup{K14n4789}\par
\textup{K14n4912}\par
\textup{K14n5071}\par
\textup{K14n5242}\par
\textup{K14n5421}\par
\textup{K14n5672}\par
\textup{K14n5928}\par
\textup{K14n6235}\par
\textup{K14n6350}\par
\textup{K14n7064}\par
\textup{K14n8816}\par
\textup{K14n9177}\par
\textup{K14n9288}\par
\textup{K14n10168}\par
\textup{K14n10618}\par
\textup{K14n10716}\par
\textup{K14n11005}\par
\textup{K14n11387}\par
\textup{K14n11539}\par
\textup{K14n11572}\par
\textup{K14n11976}\par
\textup{K14n12353}\par
\textup{K14n13160}\par
\textup{K14n13495}\par
\textup{K14n13605}\par
\textup{K14n13899}\par
\textup{K14n13909}\par
\textup{K14n14144}\par
\textup{K14n14183}\par
\textup{K14n14525}\par
\textup{K14n14615}\par
\textup{K14n14643}\par
\textup{K14n15026}\par
\textup{K14n15455}\par
\textup{K14n15642}\par
\textup{K14n15722}\par
\textup{K14n15830}\par
\textup{K14n15877}\par
\textup{K14n16318}\par
\textup{K14n16457}\par
\textup{K14n16680}\par
\textup{K14n16967}\par
\textup{K14n17032}\par
\textup{K14n17527}\par
\textup{K14n17871}\par
\textup{K14n17959}\par
\textup{K14n18129}\par
\textup{K14n18385}\par
\textup{K14n18464}\par
\textup{K14n18601}\par
\textup{K14n18796}\par
\textup{K14n18828}\par
\textup{K14n18970}\par
\textup{K14n19101}\par
\textup{K14n19485}\par
\textup{K14n19536}\par
\textup{K14n19689}\par
\textup{K14n20075}\par
\textup{K14n20766}\par
\textup{K14n21948}\par
\textup{K14n22103}\par
\textup{K14n22183}\par
\textup{K14n22269}\par
\textup{K14n22409}\par
\textup{K14n22512}\par
\textup{K14n22796}\par
\textup{K14n22949}\par
\textup{K14n23076}\par
\textup{K14n23505}\par
\textup{K14n23568}\par
\textup{K14n23986}\par
\textup{K14n24152}\par
\textup{K14n24497}\par
\textup{K14n24499}\par
\textup{K14n24527}\par
\textup{K14n24577}\par
\textup{K14n24689}\par
\textup{K14n24737}\par
\textup{K14n26226}\par
\textup{K14n26235}\par
\textup{K14n26240}\par
\textup{K14n26870}\par
\textup{K14n26875}\par
\textup{K14n26905}\par
\textup{K14n26985}\par
\textup{K14n27385}\par
\textup{K14n27399}\par
\end{multicols}
\noindent\begin{minipage}{\textwidth}
\hrule height\heavyrulewidth\vspace{3pt}
\raggedright\setlength{\parindent}{0pt}\setlength{\parskip}{0pt}
\strut For every knot listed, $h_\lambda=3$, so $\gtop\ge3$ (Corollary~\ref{cor:criteria}(c)).  The names are those of the table of Hoste and Thistlethwaite in SnapPy \cite{SnapPy}, which \cite{LS} also uses: \textup{K14a}$i$ and \textup{K14n}$i$ are the $i$th alternating and the $i$th non-alternating knot with $14$ crossings.\strut\par
\end{minipage}\par}

\Needspace{8\baselineskip}
\subsection{Knots whose value remains open}\label{app:openleft}
For the following $\nOpenLeft$ knots $\gfour\in\{1,2\}$ remains open (\S\ref{sec:values}).  So does $\gfour(\kn{11n}{176})$, for which \KI{} lists $2$ (\S\ref{sec:corrections}).  The
entries are ordered down the columns.

{\footnotesize
\setlength{\columnsep}{4pt}
\setlength{\multicolsep}{3pt}
\par\vspace{3pt}\hrule height\heavyrulewidth
\begin{multicols}{8}
\raggedright\setlength{\parindent}{0pt}\setlength{\parskip}{0pt}
$\kn{11n}{36}$\par
$\kn{11n}{40}$\par
$\kn{11n}{159}$\par
$\kn{11n}{166}$\par
$\kn{11n}{177}$\par
$\kn{11n}{178}$\par
$\kn{12a}{116}$$^{\mathrm{k}}$\par
$\kn{12a}{571}$\par
$\kn{12a}{884}$$^{\mathrm{k}}$\par
$\kn{12a}{933}$\par
$\kn{12a}{948}$\par
$\kn{12a}{976}$\par
$\kn{12a}{1263}$\par
$\kn{12n}{86}$\par
$\kn{12n}{208}$$^{\mathrm{k}}$\par
$\kn{12n}{219}$\par
$\kn{12n}{231}$\par
$\kn{12n}{255}$\par
$\kn{12n}{371}$\par
$\kn{12n}{615}$\par
$\kn{12n}{618}$\par
$\kn{12n}{730}$\par
$\kn{12n}{764}$$^{\mathrm{k}}$\par
$\kn{12n}{775}$\par
$\kn{12n}{808}$\par
$\kn{12n}{811}$\par
$\kn{13a}{27}$\par
$\kn{13a}{57}$$^{\mathrm{k}}$\par
$\kn{13a}{77}$\par
$\kn{13a}{127}$$^{\mathrm{k}}$\par
$\kn{13a}{274}$\par
$\kn{13a}{281}$\par
$\kn{13a}{313}$\par
$\kn{13a}{430}$\par
$\kn{13a}{488}$\par
$\kn{13a}{903}$\par
$\kn{13a}{919}$\par
$\kn{13a}{1351}$\par
$\kn{13a}{1668}$\par
$\kn{13a}{1835}$\par
$\kn{13a}{2028}$\par
$\kn{13a}{2469}$\par
$\kn{13a}{2562}$\par
$\kn{13a}{2575}$\par
$\kn{13a}{3050}$\par
$\kn{13a}{3051}$$^{\mathrm{k}}$\par
$\kn{13a}{3055}$\par
$\kn{13a}{3056}$\par
$\kn{13a}{3088}$\par
$\kn{13a}{3107}$$^{\mathrm{k}}$\par
$\kn{13a}{3147}$$^{\mathrm{k}}$\par
$\kn{13a}{3168}$\par
$\kn{13a}{3191}$\par
$\kn{13a}{3365}$\par
$\kn{13a}{3375}$\par
$\kn{13a}{3396}$\par
$\kn{13a}{3510}$$^{\mathrm{k}}$\par
$\kn{13a}{3517}$\par
$\kn{13a}{4000}$\par
$\kn{13a}{4010}$\par
$\kn{13a}{4031}$$^{\mathrm{k}}$\par
$\kn{13a}{4147}$\par
$\kn{13a}{4154}$$^{\mathrm{k}}$\par
$\kn{13a}{4183}$\par
$\kn{13a}{4523}$\par
$\kn{13a}{4693}$$^{\mathrm{k}}$\par
$\kn{13n}{39}$\par
$\kn{13n}{63}$\par
$\kn{13n}{205}$\par
$\kn{13n}{256}$\par
$\kn{13n}{295}$\par
$\kn{13n}{352}$\par
$\kn{13n}{453}$\par
$\kn{13n}{456}$\par
$\kn{13n}{514}$\par
$\kn{13n}{522}$\par
$\kn{13n}{533}$$^{\mathrm{k}}$\par
$\kn{13n}{541}$\par
$\kn{13n}{580}$\par
$\kn{13n}{667}$\par
$\kn{13n}{736}$\par
$\kn{13n}{764}$\par
$\kn{13n}{804}$\par
$\kn{13n}{845}$\par
$\kn{13n}{849}$\par
$\kn{13n}{864}$\par
$\kn{13n}{869}$$^{\mathrm{k}}$\par
$\kn{13n}{903}$\par
$\kn{13n}{909}$\par
$\kn{13n}{938}$$^{\mathrm{k}}$\par
$\kn{13n}{1069}$\par
$\kn{13n}{1094}$\par
$\kn{13n}{1213}$\par
$\kn{13n}{1242}$\par
$\kn{13n}{1244}$\par
$\kn{13n}{1255}$$^{\mathrm{k}}$\par
$\kn{13n}{1258}$\par
$\kn{13n}{1263}$$^{\mathrm{k}}$\par
$\kn{13n}{1350}$$^{\mathrm{k}}$\par
$\kn{13n}{1389}$\par
$\kn{13n}{1529}$$^{\mathrm{k}}$\par
$\kn{13n}{1686}$\par
$\kn{13n}{1720}$\par
$\kn{13n}{1736}$\par
$\kn{13n}{1781}$\par
$\kn{13n}{1782}$$^{\mathrm{k}}$\par
$\kn{13n}{1822}$\par
$\kn{13n}{1870}$\par
$\kn{13n}{1908}$\par
$\kn{13n}{1909}$$^{\mathrm{k}}$\par
$\kn{13n}{1949}$$^{\mathrm{k}}$\par
$\kn{13n}{1965}$\par
$\kn{13n}{2083}$\par
$\kn{13n}{2184}$\par
$\kn{13n}{2186}$\par
$\kn{13n}{2237}$\par
$\kn{13n}{2413}$$^{\mathrm{k}}$\par
$\kn{13n}{2425}$\par
$\kn{13n}{2575}$$^{\mathrm{k}}$\par
$\kn{13n}{2658}$$^{\mathrm{k}}$\par
$\kn{13n}{2743}$\par
$\kn{13n}{2829}$\par
$\kn{13n}{2855}$\par
$\kn{13n}{2979}$\par
$\kn{13n}{2980}$\par
$\kn{13n}{3045}$$^{\mathrm{k}}$\par
$\kn{13n}{3080}$\par
$\kn{13n}{3117}$\par
$\kn{13n}{3299}$\par
$\kn{13n}{3307}$\par
$\kn{13n}{3352}$\par
$\kn{13n}{3357}$$^{\mathrm{k}}$\par
$\kn{13n}{3358}$\par
$\kn{13n}{3367}$\par
$\kn{13n}{3382}$\par
$\kn{13n}{3399}$\par
$\kn{13n}{3445}$\par
$\kn{13n}{3449}$\par
$\kn{13n}{3455}$$^{\mathrm{k}}$\par
$\kn{13n}{3479}$\par
$\kn{13n}{3496}$\par
$\kn{13n}{3505}$\par
$\kn{13n}{3534}$\par
$\kn{13n}{3548}$\par
$\kn{13n}{3628}$\par
$\kn{13n}{3662}$\par
$\kn{13n}{3663}$\par
$\kn{13n}{3664}$\par
$\kn{13n}{3704}$\par
$\kn{13n}{3804}$\par
$\kn{13n}{3840}$$^{\mathrm{k}}$\par
$\kn{13n}{3848}$\par
$\kn{13n}{3864}$\par
$\kn{13n}{3893}$\par
$\kn{13n}{3898}$\par
$\kn{13n}{3941}$$^{\mathrm{k}}$\par
$\kn{13n}{3944}$\par
$\kn{13n}{3951}$\par
$\kn{13n}{3958}$\par
$\kn{13n}{3981}$\par
$\kn{13n}{4000}$\par
$\kn{13n}{4035}$$^{\mathrm{k}}$\par
$\kn{13n}{4062}$\par
$\kn{13n}{4091}$\par
$\kn{13n}{4093}$\par
$\kn{13n}{4104}$\par
$\kn{13n}{4117}$\par
$\kn{13n}{4134}$\par
$\kn{13n}{4148}$\par
$\kn{13n}{4177}$\par
$\kn{13n}{4192}$\par
$\kn{13n}{4205}$$^{\mathrm{k}}$\par
$\kn{13n}{4238}$\par
$\kn{13n}{4260}$$^{\mathrm{k}}$\par
$\kn{13n}{4283}$\par
$\kn{13n}{4290}$\par
$\kn{13n}{4323}$\par
$\kn{13n}{4349}$\par
$\kn{13n}{4358}$\par
$\kn{13n}{4371}$$^{\mathrm{k}}$\par
$\kn{13n}{4385}$$^{\mathrm{k}}$\par
$\kn{13n}{4387}$\par
$\kn{13n}{4389}$$^{\mathrm{k}}$\par
$\kn{13n}{4391}$\par
$\kn{13n}{4450}$\par
$\kn{13n}{4452}$$^{\mathrm{k}}$\par
$\kn{13n}{4454}$\par
$\kn{13n}{4467}$\par
$\kn{13n}{4468}$\par
$\kn{13n}{4478}$\par
$\kn{13n}{4495}$$^{\mathrm{k}}$\par
$\kn{13n}{4553}$\par
$\kn{13n}{4554}$$^{\mathrm{k}}$\par
$\kn{13n}{4572}$\par
$\kn{13n}{4577}$\par
$\kn{13n}{4600}$\par
$\kn{13n}{4663}$\par
$\kn{13n}{4688}$\par
$\kn{13n}{4708}$$^{\mathrm{k}}$\par
$\kn{13n}{4726}$$^{\mathrm{k}}$\par
$\kn{13n}{4733}$\par
$\kn{13n}{4750}$\par
$\kn{13n}{4756}$$^{\mathrm{k}}$\par
$\kn{13n}{4773}$\par
$\kn{13n}{4775}$$^{\mathrm{k}}$\par
$\kn{13n}{4778}$$^{\mathrm{k}}$\par
$\kn{13n}{4799}$\par
$\kn{13n}{4818}$$^{\mathrm{k}}$\par
$\kn{13n}{4828}$\par
$\kn{13n}{4848}$\par
$\kn{13n}{4852}$$^{\mathrm{k}}$\par
$\kn{13n}{4860}$$^{\mathrm{k}}$\par
$\kn{13n}{4886}$\par
$\kn{13n}{4888}$$^{\mathrm{k}}$\par
$\kn{13n}{4889}$$^{\mathrm{k}}$\par
$\kn{13n}{4904}$\par
$\kn{13n}{4916}$$^{\mathrm{k}}$\par
$\kn{13n}{4926}$$^{\mathrm{k}}$\par
$\kn{13n}{4935}$\par
$\kn{13n}{4999}$$^{\mathrm{k}}$\par
$\kn{13n}{5025}$$^{\mathrm{k}}$\par
$\kn{13n}{5037}$$^{\mathrm{k}}$\par
$\kn{13n}{5040}$\par
$\kn{13n}{5047}$\par
$\kn{13n}{5079}$\par
\end{multicols}
\noindent\begin{minipage}{\textwidth}
\hrule height\heavyrulewidth\vspace{3pt}
\raggedright\setlength{\parindent}{0pt}\setlength{\parskip}{0pt}
\strut Unmarked: \KI{} range $[1,2]$.\strut\par
\strut $^{\mathrm{k}}$ \KI{} range $[1,3]$; the knot bounds a punctured Klein bottle made of band moves (Proposition~\ref{prop:upper}(b)), so $\gfour\in\{1,2\}$.\strut\par
\end{minipage}\par}

\end{document}